\documentclass[a4paper,11pt,reqno]{amsart}
\usepackage{amsmath,amssymb,amsthm} %Some extra symbols
\usepackage{graphics} %If you want to include postscript graphics
\usepackage{enumerate}

\usepackage[colorlinks,citecolor=blue]{hyperref}
\usepackage{tikz}

\usetikzlibrary{patterns, positioning}
\usepackage{epsfig}
\usepackage{color}

\numberwithin{equation}{section}

\newtheorem{theorem}{Theorem}[section]
\newtheorem{lemma}[theorem]{Lemma}
\newtheorem{proposition}[theorem]{Proposition}

\newtheorem{hypothesis}[theorem]{Hypothesis}

\theoremstyle{definition}

\newtheorem{remark}[theorem]{Remark}

\newcommand{\hel} {
\hskip2.5pt{\vrule height7pt width.5pt depth0pt}
\hskip-.2pt\vbox{\hrule height.5pt width7pt depth0pt}
\, }
\newcommand{\restr}{\hel}
\newcommand{\R}{\mathbb{R}}
\newcommand{\N}{\mathbb{N}}

\newcommand{\Z}{\mathbb{Z}}
\newcommand{\vhi}{\varphi}
\newcommand{\eps}{\varepsilon}

\newcommand{\be}{\begin{equation}}
\newcommand{\ee}{\end{equation}}

\newcommand\lt{\left}
\newcommand\rt{\right}

\def\HH{\mathcal{H}}
\def\les{\lesssim}
\def\ges{\gtrsim}

\def\M{\mathcal{M}}

\def\reps{\widehat{r}}
\def\m{m}

\def\X{\mathcal{X}}

\def\bS{\mathbb{S}}
\def\cP{\mathcal{P}}

\newcommand{\ang}[1]{\lt\langle #1 \rt\rangle}
\newcommand{\abs}[1]{\left| #1 \right|}

\renewcommand{\tilde}{\widetilde}
\def\diam{{\rm diam}}
\def\spt{ {\rm Spt}\,}
\title{Quantitative rigidity of the Wasserstein contraction under convolution}

\author{Max Fathi}
\address{Université Paris Cité and Sorbonne Université, CNRS, Laboratoire Jacques-Louis Lions and Laboratoire de Probabilit\'es, Statistique et Mod\'elisation, F-75013 Paris, France\\
\newline and DMA, École normale supérieure, Université PSL, CNRS, 75005 Paris, France \\
\newline and Institut Universitaire de France }
\email{mfathi@lpsm.paris}

\author{Michael Goldman}
\address{CMAP, CNRS, \'Ecole polytechnique, Institut Polytechnique de Paris, 91120 Palaiseau,
France}
\email{michael.goldman@cnrs.fr}
\date{}

\author{Daniel Tsodyks}
\address{Department of Mathematics, Weizmann Institute of Science, Rehovot 7610001, Israel}
\email{danieltsod@gmail.com}

\begin{document}

\begin{abstract}
 The aim of this paper is to investigate the contraction properties of $p$-Wasserstein distances with respect to convolution  in Euclidean spaces both qualitatively and quantitatively. We connect this question to the question of uniform convexity of the Kantorovich functional on which there was substantial recent progress (mostly for $p=2$ and partially for $p>1$). Motivated by this connection we extend these uniform convexity results to the case $p=1$, which is of independent interest.
\end{abstract}
\maketitle

\section{Introduction}
Let $\lambda,\mu$ be probability measures on $\R^n$ and let $\rho_\eps : \R ^ n \rightarrow \R$ be a standard convolution kernel; that is, $\rho_\eps(x)=\eps^{-n}\rho(x/\eps)$ for some function $\rho\ge 0$ satisfying $\int_{\R^n} \rho dx=1$. Denote $\lambda_\eps=\rho_\eps\ast \lambda$ and $\mu_\eps=\rho_\eps\ast \mu$. It is well-known, see e.g. \cite[Lemma 5.2]{Santam} as well as the proof of Theorem \ref{theo:equal:intro} below, that for every $p\ge 1$, and every $\eps>0$,
\begin{equation}\label{eq:contraction}
 W_p(\lambda_\eps, \mu_\eps) \le W_p(\lambda, \mu).
\end{equation}
Here $W_p$ stands for the standard $p-$Wasserstein distance. For applications of this inequality as well as variants in more general geometries, we refer for instance to  \cite{von2005transport,AmStTr16,cosenza2024new,chen2022asymptotics,bolley2014dimensional}. {Two classical situations where this inequality is relevant are
\begin{itemize}
\item Contractive estimates for Wasserstein distances serve as a good definition of non-negative Ricci curvature in various settings \cite{von2005transport, ollivier2009}, due to the fact that decay of Wasserstein distances along the heat flow is a characterisation of nonnegative curvature (and the heat flow is a convolution kernel in the Euclidean setting);

\item Often in statistics one wishes to estimate a Wasserstein distance $W_p(\lambda, \mu)$ to compare an empirical distribution to some theoretical one, but only has access to samples from noisy distributions  $\lambda \ast \rho_\eps$, where $\rho$ is the distribution of the noise, and $\eps$ a small size parameter. Inequality \eqref{eq:contraction} implies that independent noise makes us underestimate the distance.
\end{itemize}}

Setting
\begin{equation}\label{def:deltaepsmunu}
 \delta_\eps(\lambda,\mu)=W^p_p(\lambda,\mu)- W^p_p( \lambda_\eps, \mu_\eps),
\end{equation}
the contraction property \eqref{eq:contraction} may be rephrased as $\delta_\eps(\lambda,\mu)\ge 0$. The first question we aim at understanding is the equality case $\delta_\eps(\lambda,\mu)=0$. Unless specified otherwise, we will assume throughout the paper that $\lambda$ (but not necessarily $\mu$) is absolutely continuous with respect to the Lebesgue measure; we will also assume that $\inf_{B_1} \rho > 0$ (we denote by $B_R(x)$ the open ball of radius $R>0$ centered at $x\in \R^n$ and $B_R=B_R(0)$). Notice that if $\spt \rho$ is compact and the distances between the connected components of $\X= \spt \lambda$ are bounded away from zero, then for $\eps$ small enough the convolution acts independently on each connected component of $\X$. Therefore, we will most often assume that $\X$ is connected. For a measure $\mu$ and $z\in \R^n$, we define the translated measure $\mu^z$ by $\mu^z(A)=\mu(A-z)$ for $A\subset \R^n$.
\begin{theorem}\label{theo:equal:intro}
Let $\lambda$ and $\mu$ be probability measures on $\R^n$. Assume that either $\spt \rho=\R^n$ or $\X$ is connected and moreover that there exist optimal Kantorovich potentials for $W_p(\lambda,\mu)$ and $W_p(\lambda_\eps,\mu_\eps)$. Then  $\delta_\eps(\lambda,\mu)=0$ if and only if:
\begin{itemize}
\item For $p>1$, there exists $z\in \R^n$ such that $\mu=\lambda^z$.
\item For $p=1$,  there exists $e\in \bS^{n-1}$ such that whenever $\varphi\in L^1(\mu)\cap L^1(\lambda)$ is monotone in the direction $e$, i.e. $\varphi(x+te)\ge \varphi(x)$ for all $x\in \R^n$ and $t\ge 0$, we have
\begin{equation}\label{condmonotonephi}
 \int_{\R^n} \varphi d\mu\ge \int_{\R^n} \varphi d\lambda.
\end{equation}
\end{itemize}
\end{theorem}
{
\begin{remark}
The fact that the case $p=1$ is different can easily be understood in dimension one. If we consider two measures with disjoint supports $[a,b]$ and $[c,d]$ with $b < c$, then if $\rho$ is compactly supported, for any $\eps$ small enough, $\lambda_\eps$ and $\mu_{\eps}$ also have disjoint supports, and for any coupling $\pi_\eps$ of the two perturbed measures
$$W_1(\lambda_\eps, \mu_\eps) = \int_{\R\times\R}{|x-y|d\pi_{\eps}} = \left|\int_{\R\times\R}{(x-y) d\pi_\eps}\right| = \left|\int_\R{xd\lambda} - \int_\R{yd\mu}\right|=W_1(\lambda,\mu).$$
\end{remark}}
\begin{remark}\label{rem:equiv}
In both cases the conclusion may be seen to be equivalent to the fact that the gradient of the optimal Kantorovich potential for $W_p(\lambda,\mu)$ is constant, see the proof of Theorem \ref{theo:equal:intro} in Section \ref{sec:convol}. Moreover, when $p=1$,
letting $p_e$ be the orthogonal projection on $e^\perp$, \eqref{condmonotonephi} is equivalent to the following property. We have   $p_e\#\lambda=p_e\#\mu=\eta$ for some measure $\eta$ and, disintegrating $\lambda=\lambda_y \otimes \eta$ and $\mu= \mu_y\otimes \eta$, $\mu_y$ is stochastically dominated by $\lambda_y$ on $p_e^{-1}(y)$ for $\eta-$almost every $y$. To see that \eqref{condmonotonephi} implies $p_e\#\lambda=p_e\#\mu$, it is enough to consider all the  test functions $\varphi$ which do not depend on $\ang{x,e}$, since  for these we have equality in \eqref{condmonotonephi} (by testing the inequality with  both $\varphi$ and $-\varphi$). The second property then follows. Another informal way of interpreting \eqref{condmonotonephi} is that transport between $\lambda$ and $\mu$ happens only in the direction $e$.
\end{remark}

{
\begin{remark}
In a statistical context, one may think of this rigidity result as stating that when $p > 1$, the only situation in which there is no loss from estimating the Wasserstein distance on noisy data rather than directly is when considering a translation model.
\end{remark}
}
% \begin{theorem}\label{theo:equal:intro}
% Let $\lambda$ and $\mu$ be probability measures on $\R^n$. Assume that either $\spt \rho=\R^n$ or $\X$ is connected and moreover that there exists optimal Kantorovich potentials for $W_p(\lambda,\mu)$ and $W_p(\lambda_\eps,\mu_\eps)$. Then  $\delta_\eps(\lambda,\mu)=0$ if and only if:
% \begin{itemize}
% \item For $p>1$, there exists $z\in \R^n$ such that $\mu=\lambda^z$.
% \item For $p=1$,  there exists $e\in \bS^{n-1}$ such that letting $p_e$ be the orthogonal projection on $e^\perp$, $p_e\#\lambda=p_e\#\mu=\eta$. Moreover, for the disintegration $\lambda=\lambda_y \otimes \eta$ and $\mu= \mu_y\otimes \eta$, $\lambda_y$ is stochastically dominated by $\mu_y$ on $p_e^{-1}(y)$ for $\eta-$almost every $y$.
% \end{itemize}
% \end{theorem}
We then consider the quantitative counterpart of this question. More precisely, we ask the following question: for $p>1$, if $\delta_\eps(\lambda,\mu)$ is small, does it imply that $\mu$ is almost a translate of $\lambda$? Similarly, if $p=1$ and $\delta_\eps(\lambda,\mu)$ is small, does it mean that there is a direction $e\in \bS^{n-1}$ such that \eqref{condmonotonephi} almost holds?\\
Currently, the answer to these questions depends on the hypothesis which are made on $\lambda$ and $\mu$. Rather than giving a general statement, let us give an example of what we can achieve.
\begin{theorem}\label{theo:introsimplequadratic}
 Let $p=2$, and assume that $\X=\spt \lambda$ is  bounded, Lipschitz and connected. If $m\le \lambda\le M$ on $\X$ for some $0<m<M$ and $R, \eps>0$ is such that $\spt\lambda_\eps\cup\spt \mu_\eps \subset B_R$, then there exists $C=C(\rho,R,\X,m,M)>0$ such that
 \begin{equation}\label{statement:simplequadratic}
  \min_{z\in \R^n} W_2^2(\lambda,\mu^z)\le C \eps^{-(n+1)} \delta_\eps^{\frac{1}{3}}(\lambda,\mu).
 \end{equation}
\end{theorem}
\begin{remark}\label{remrho}
 Let us make an important observation regarding the convolution kernel $\rho$. For most of our results it will be important that the measures we consider both before and after convolution  have bounded support. In that case we will thus restrict ourselves to convolution kernels which have themselves bounded support (as in the theorem above). To illustrate  the applicability of the method beyond the compact case, we will also consider in Theorem \ref{theo:stabgauss} the case when $\lambda$ (but not $\mu$) is Gaussian and $\rho_\eps$ is the heat kernel.
\end{remark}
\begin{remark}\label{remoptim}
 Since $\delta_\eps(\lambda,\mu)\le W_p^p(\lambda,\mu)$, we see that the best possible exponent we can hope for on the right-hand side of \eqref{statement:simplequadratic} is one. As seen in Theorem \ref{theo:conditional}, we would obtain this sharp exponent provided the strong convexity estimates from \cite{DelalMeri} could be  improved accordingly. Regarding the dependence of the right-hand side of \eqref{statement:simplequadratic} in $\eps$ we do not claim any optimality, see in particular Remark \ref{rem:gaussian case}. In this respect notice also that by \eqref{eq:tauplambdaConv} of Lemma \ref{lem:tauplambda} we get an improved dependence on $\eps$ if $\X$ is convex.
\end{remark}

We now discuss the proof of Theorem \ref{theo:introsimplequadratic}.
Our first main insight is that this question is closely related to the question of stability in optimal transportation. Indeed, from our proof of Theorem \ref{theo:equal:intro}, we deduce that if $\delta_\eps(\lambda,\mu)\ll1$, then the optimal Kantorovich potentials for $ W_2( \lambda_\eps, \mu_\eps)$ are almost optimal potentials for the dual formulation of  $W_2(\lambda,\mu)$, and similarly for all the translations of $\lambda$ and $\mu$. Thanks to recent progress on the strong convexity properties of the Kantorovich functional in the case of the quadratic cost, see \cite{gigli2011holder,hutter2021minimax,DelalMeri,merigot2020quantitative,LetMer,FordThesis}, we can deduce that the optimal transport map $T_\eps$ for $W_2(\lambda_\eps,\mu_\eps)$ is itself close to the optimal transport map $T$ for $W_2(\lambda,\mu)$. We thus obtain a control on the quantity (which we define here for all $p\ge 1$)
\begin{equation}\label{defLambda}
  \Lambda_\eps(\xi,f)=\int_{\R^n\times \R^n} |\xi(y)-f(x)|^p\rho_\eps(x-y)dx d\lambda(y),
\end{equation}
where $\xi(y)=T(y)-y$ and $f(x)=T_\eps(x)-x$. Our second main insight is that if $\Lambda_\eps$ is small then $\xi$ is close to being constant. For the following theorem, We let $\tau_{p,\lambda}(\eps)$ be the quantity defined in \eqref{deftaup}, see also Lemma \ref{lem:tauplambda}.
\begin{theorem}\label{twopointcompact}
 Let $\lambda$ be such that Hypothesis \ref{hyp:lambda} holds for $r\leq\eps$. For $\xi\in L^1(\lambda)$, letting,
\[z=\int_{\R^n} \xi d\lambda, \]
we have for every $f$ and every $p\ge 1$,
\begin{equation}\label{eq:twopointcompact}
 \int_{\R^n}|\xi-z|^p d\lambda\les M_0(r)(\eps/r)^n\tau_{p,\lambda}(r) \Lambda_\eps(\xi,f).
\end{equation}
Here the implicit constant depends only on the constant $\eta$ from Hypothesis \ref{hyp:lambda} and on $\inf_{B_1} \rho$.
\end{theorem}
In order to discuss the proof of \eqref{eq:twopointcompact}, let us first consider the case $\Lambda_\eps(\xi,f)=0$. In that case, for every $x\in \R^n$ and every $y\in \spt \lambda \cap  B_\eps(x)$ we have $\xi(y)=f(x)$. In particular, if for some $x,x'$, we have $\spt \lambda \cap  B_\eps(x)\cap B_\eps(x')\neq \emptyset$, then  $\xi(y)=f(x)=f(x')$ in $\spt \lambda \cap (B_\eps(x)\cup B_\eps(x'))$. Covering $\spt \lambda$ by overlapping balls of radius $\eps$, we find that $\xi$ is constant on $\spt \lambda$ (the assumption $\tau_{p,\lambda}(\eps)<\infty$ will turn out to imply some $\eps-$connectedness of $\spt \lambda$). From this argument, we learn that $\Lambda_\eps(\xi,f)$ controls how far is $\xi$ from being constant in balls of radius $\eps$. In order to prove that $\xi$ is globally close to a  constant we use a gluing argument which is reminiscent of some proofs of the Poincar\'e inequality, see e.g. \cite{grigor2005stability,hajlasz2000sobolev}. This is also very similar in spirit to the proofs of quantitative stability of optimal transport maps in \cite{kitagawa2025stability, LetMer}. Let us also point out that the quantity $\Lambda_\eps(\xi,f)$ is closely related to non-local characterizations of Sobolev spaces by Bourgain-Brézis-Mironescu, see \cite{Bourgain01,brezis2002recognize} as well as Remark \ref{rmk:BBM} below. In particular Theorem \ref{twopointcompact} generalizes the main result of \cite{Ponce04} with a very different proof which is somewhat closer to arguments from \cite{goldman2023non}.

Let us now comment on extensions of Theorem \ref{theo:introsimplequadratic}. As explained above, besides Theorem \ref{twopointcompact}, one of the main building blocks in our proof of Theorem \ref{theo:introsimplequadratic} is the uniform convexity of the Kantorovich functional. Despite recent progress on this question, since at the moment of writing this paper sharp results do not seem to be available, we rather write our main theorem as a conditional result as in \cite{carlier2024quantitative}, see Theorem \ref{theo:conditional}. Let us point out that while many results are available regarding strong monotonicity of the gradient of the Kantorovich functional, see \cite{ gigli2011holder,DelalMeri,merigot2020quantitative,LetMer,kitagawa2025stability} for $p=2$ and \cite{MisTre} for $p>1$, we actually need as in \cite{carlier2024quantitative,HuGoTre} the (seemingly stronger) strong convexity of the Kantorovich functional itself. This in turn is currently known only under either strong regularity properties of the potentials, see \cite{gigli2011holder,hutter2021minimax,merigot2020quantitative} or in the setting of Theorem \ref{theo:introsimplequadratic}.  We finally point out that since $\tau_{p,\lambda}(\eps)<\infty$ morally implies that $\spt \lambda$ is bounded, see Remark \ref{rem:taubounded}, in order to show that our strategy of proof goes beyond the compact case, we consider in Theorem \ref{theo:stabgauss} also the case when $\lambda$ is a Gaussian measure.\\

Our second main contribution of the paper is a quantitative stability result of Kantorovich potentials when $p=1$. Since in this case, we have $W_1(\lambda+\nu,\mu+\nu)=W_1(\lambda,\mu)$ whenever $\lambda(\R^n)=\mu(\R^n)$, we will always assume that $\lambda\perp \mu$. We will actually make the stronger hypothesis that $\lambda$ and $\mu$ have disjoint support, see Remark \ref{rem:support} for a discussion of this assumption.
\begin{theorem}
\label{thm1}
    Let $\lambda, \mu$ be probability measures with disjoint supports both contained in $ B_R$ for some $R>0$, and let $\X=\spt \lambda$. Assume $\partial \X$ is rectifiable and $m \leq \lambda \leq M$ on $\X$  for some  $0 < m, M < \infty$. Let $\psi$ be the Kantorovich potential associated with $W_1(\lambda, \mu)$, and $\phi$ be any $1$-Lipschitz function. Then, for every $\alpha>3$, there exist a constant $C=C(\alpha, R,\X,m,M)>0$, such that
    \begin{equation}\label{eq:thm1}
        \|\nabla \psi - \nabla \phi\|_{L ^ 1(\lambda)}^\alpha \le C \int_{\R^n} (\psi - \phi)d(\mu - \lambda) =C \left(W_1(\lambda,\mu)-\int_{\R^n}  \phi d(\mu - \lambda) \right).
    \end{equation}
\end{theorem}
We recall that a compact set $K\subset \R^n$ is called rectifiable if there exists a Lipschitz function $f:\R^{n-1}\to \R^n$ and a compact set $K_0$ such that $K=f(K_0)$. Let us point out that while there is no reason to believe that the condition $\alpha>3$ is sharp, \eqref{eq:thm1} cannot hold for $\alpha <2$, see Remark \ref{optimality}.\\
The proof of Theorem \ref{thm1} is very different from the proof from \cite{DelalMeri} in the quadratic case. The starting point is Lemma \ref{lem:centralstab1} where a stability estimate of the form \eqref{eq:thm1} is obtained via a relatively simple argument but at the price of replacing the measure $\lambda$ on the left-hand side of \eqref{eq:thm1} by the transport density $\sigma$ (see \eqref{eq:transportdens} for the definition). The difficulty is then to estimate $\sigma$ from below by $\lambda$, see in particular Lemma \ref{lem1}. Let us point out that in Theorem \ref{thm2} we extend the validity of \eqref{eq:thm1} to  some measures which are not bounded from below and might have unbounded support.\\

 Equipped with Theorem \ref{thm1}, we can prove an analog of Theorem \ref{theo:introsimplequadratic} for $p=1$.
\begin{theorem}\label{theo:introsimplelinear}
 Let $\lambda, \mu$ be probability measures with disjoint supports both contained in $ B_R$ for some $R>0$ and let $\X=\spt \lambda$. Assume that $ \X$ is Lipschitz and  $m \leq \lambda \leq M$ on $\X$  for some  $0 < m, M < \infty$. For every $\alpha>3$, there exist $C=C(\alpha,\X,m, M,R)>0$ and  $e\in \bS^{n-1}$ such that for every  $1-$Lipschitz and  monotone  function $\vhi$ in the direction $e$,
 \begin{equation}\label{quantmonot}
  \int_{\R^n} \vhi d\lambda- \int_{\R^n} \vhi d\mu\le C\eps^{-n} \delta_\eps^{\frac{1}{\alpha}}(\lambda,\mu).
 \end{equation}
As a consequence (recall the notation from Remark \ref{rem:equiv})
\begin{equation}\label{stabmarg}
 W_1(p_e\#\lambda,p_e\#\mu)\le C \eps^{-n} \delta_\eps^{\frac{1}{\alpha}}(\lambda,\mu).
\end{equation}

\end{theorem}

The paper is organized as follows. In Section \ref{sec:notation} we gather the notation we use throughout the paper and recall some facts about optimal transportation. In Section \ref{sec:2point}, we prove Theorem \ref{twopointcompact}. In Section \ref{sec:stablin} we investigate the stability properties of Kantorovich potentials for $p=1$. In particular we prove there Theorem \ref{thm1}. Finally, in Section \ref{sec:convol} we turn to our main objective which is the study of \eqref{eq:contraction} and its rigidity. We first prove a qualitative version, namely Theorem \ref{theo:equal:intro} and then the quantitative counterparts for $p>1$ in Section \ref{subsec:stabp} and $p=1$ in Section \ref{subsec:stapone}.

\section{Notation and preliminaries}\label{sec:notation}
Given a set $A$, we write $\chi_A$ for its indicator function. We write $\ang{x, y}$ for the standard scalar product in $\R^n$ and $\abs{x} = \sqrt{ \ang{x, x}}$ for the Euclidean norm. We always endow $\R^n$ with the Euclidean distance. The notation $A\les B$ means that there exists a constant $C>0$, such that $A\le C B$, where unless specified otherwise, $C$ depends on the dimension $n$ and $p$. We write $A \les_q B$ to indicate an additional dependence on the parameter, or set of parameters, $q$. We write $A\approx B$ if both $A\les B$ and $B\les A$. We use the notation $A\ll B$ as an hypothesis; it means that there exists some (typically small) $\eps>0$ such that if $A\le \eps B$, then the conclusion holds. For $p>1$, we write $p'$ for the H\"older conjugate of $p$, which is the unique $p' > 1$ for which $1/p+ 1/p'=1$. For $e\in \bS^{n-1}$, we set $p_e(x)=x- \ang{x, e}e$ to be the orthogonal projection on $e^\perp$. For a closed set $\X \subseteq \R ^ n$, we write $\cP(\X)$ for the set of probability measures supported on $\X$, and $\M(\X)$ for the set of positive Radon measures supported on $\X$. If a measure $\lambda$ is absolutely continuous with respect to the Lebesgue measure, by abuse of notation, we will identify it with its density, so $\lambda(x)$ is the density of $\lambda$ at point $x$. Similarly, if $\lambda$ is radially symmetric we identify $\lambda(x)$ and $\lambda(|x|)$. We say that a measure $\mu$ is log-concave if $\mu=\exp(-V)$ for some convex function $V: \R^n\to \R\cup\{+\infty\}$. {Note that this is a bit more restrictive than the usual definition of log-concavity, which allows for more degenerate measures, that are supported on a hyperplane, but we shall not consider degenerately log-concave measures here.} If, moreover, $V$ is $\kappa-$strongly convex, we say that $\mu$ is $\kappa-$uniformly log concave.
\subsection{Optimal transport}
In this section we recall some basic facts regarding optimal transport, see e.g. \cite{Santam,Viltop,AmbPra}. For $p\ge 1$ and measures $\lambda,\mu$ with finite $p-$moment and $\lambda(\R^n)=\mu(\R^n)$ we set
\[
 \Pi(\lambda,\mu)=\{\pi\in \M(\R^n\times \R^n) \ : \ \pi_1=\lambda, \ \pi_2=\mu\},
\]
where $\pi_1, \pi_2$ are the marginals of $\pi$. Then, we define
\[
 W_p^p(\lambda,\mu)=\inf_{\Pi(\lambda,\mu)} \int_{\R^n\times \R^n} \frac{1}{p}|x-y|^p d\pi.
\]
If $\lambda$ is absolutely continuous, then there exists an optimal transport plan $\pi$ such that $\pi=(Id\times T)\#\lambda$ for some optimal transport map $T$ (see e.g. \cite[Theorem 7.1]{AmbPra} for the existence in the  case $p=1$). When $p>1$, $T$ and $\pi$ are unique; when $p=1$, however, there is no uniqueness, even in the one-dimensional case. Nevertheless, if $p=1$ there is a unique ray-monotone optimal transport map, see \cite[Theorem 6.4 \& Theorem 6.5]{AmbPra} as well as \cite[Theorem 3.18]{Santam} where the case of measures with compact support is covered.\\
For $t\in[0,1]$ and taking $T$ to be a an optimal transport map, we set
\begin{equation}\label{defTt}
 T_t(x)=(1-t)x +tT(x),
\end{equation}
and then
\begin{equation}\label{defmuj}
 \mu_t=T_t\#\lambda, \qquad j_t=T_t\#\lt( (T-x)\lambda\rt),
\end{equation}
so that in the distributional sense,
\[
 \partial_t \mu_t + \nabla \cdot j_t=0.
\]
By integration we get that $\bar j=\int_0^1 j_t dt$ satisfies
\[
 \nabla\cdot \bar j= \mu_0 - \mu_1 = \lambda-\mu.
\]
For $\psi:\R^n\to \R\cup\{-\infty\}$, with $\psi$ not identically equal to $-\infty$, we set
\[
 \psi^c(x)=\inf_{y\in \R^n} \frac{1}{p}|x-y|^p-\psi(y).
\]
If $\phi=\psi^c$ for some function $\psi$, we say that $\phi$ is $c-$concave. We now define the Kantorovich functional as
\begin{equation}\label{def:Kantorovichfunct}
 F_{\lambda,\mu}(\psi)=\int_{\R^n} \psi^c d\lambda +\int_{\R^n} \psi d\mu.
\end{equation}
Kantorovich duality asserts that we have
\be
\label{def:Wpdual}
 W_p^p(\lambda,\mu)=\sup_\psi F_{\lambda,\mu}(\psi).
\ee
If now $R>0$ is such that $\spt \lambda\cup \spt \mu\subset B_R$, the supremum in \eqref{def:Wpdual} is achieved by some $\psi$. Since the value of $\psi$ outside $B_R$ plays no role, we may assume that $\psi=-\infty$ outside $B_R$, in which case $\psi^c$ is Lipschitz continuous with
\begin{equation}\label{eq:Lippsic}
 \sup_{B_R} |\nabla \psi^c|\les R^{p-1}.
\end{equation}
When $p=1$, it is not hard to see that a function $\psi$ is $c$-concave if and only if it is $1$-Lipschitz. Moreover, in that case we have $\psi^c=-\psi$, so that
\begin{equation}\label{defW1dual}
 W_1(\lambda,\mu)=\sup_{\|\nabla \psi\|_\infty\le 1} \int_{\R^n} \psi d(\mu-\lambda).
\end{equation}
We now connect optimal potentials and optimal transport maps. For $p>1$, let
\begin{equation}\label{def:Phip}
 \Phi_p(z)=|z|^{p'-2} z.
\end{equation}
We then have that if $T$ is the optimal transport map and $\psi$ is an optimal Kantorovich potential,
\begin{equation}\label{PhipT}
 x-T(x)=\Phi_p(\nabla \psi^c(x)).
\end{equation}
Since $\Phi_p$ is invertible, this implies that $\nabla \psi^c$ is determined uniquely on $\X=\spt \lambda$. For $p=1$, we have that
\[
 |T(x)-x|=\psi(T(x))-\psi(x),
\]
and because $\psi$ is $1$-Lipschitz, if $T(x)\neq x$ we have
\begin{equation}\label{nabpsi}
 \nabla \psi(x)= \frac{T(x)-x}{|T(x)-(x)|},
\end{equation}
so that the mass is transported through rays of direction given by $\nabla \psi$. Moreover $\psi$ is affine (with constant gradient) on the segment  $[x,T(x)]$. We thus have
\begin{equation}\label{eq:transportdens}
 \sigma=|\bar j|=\int_0^1|j_t| dt.
\end{equation}
Notice that \eqref{nabpsi} and the fact that transport happens along rays implies that $\nabla \psi$ is determined uniquely on $\spt \sigma$.\\
\begin{lemma}\label{lem:qualitabscont}
 Let $p=1$. If  $\lambda\perp \mu$ then $\spt \lambda\subset \spt \sigma$.
\end{lemma}
\begin{proof}
 We start by noticing that since $\lambda\perp \mu$, for $\lambda$ a.e. $x$ we have $|T(x)-x|>0$. Without loss of generality, assume that $0 \in \spt \lambda$, and suppose for the sake of contradiction that $\sigma(B_r)=0$ for some $r>0$. We then have by \eqref{eq:transportdens}, and the definition \eqref{defmuj} of $j_t$,
 \[
  \sigma(B_r)=\int_0^1\int_{T_{t}^{-1}(B_r)}|T(x) - x| d\lambda(x) dt=0,
 \]
and thus for almost every $t\in(0,1)$ and $\lambda$-almost every $x\in T_{t}^{-1}(B_r)$, we get $|T(x)-x|=0$. In order to obtain a contradiction with the initial observation, let us construct a set $A\subset B_{r/2}\cap T_{t}^{-1}(B_r)$ with $\lambda(A)>0$. For $M>0$ let
\[
 A_M=\{x \ : \ |T(x)-x|\le M\}.
\]
Since $W_1(\lambda,\mu)<\infty$, by the Markov inequality we have $\lambda(A_M^c)\le \lambda(B_{r/2})/2$ provided $M$ is large enough. Setting $A=A_M\cap B_{r/2}$ we thus have
\[
 \lambda(A)\ge \lambda(B_{r/2})/2>0.
\]
Moreover, for every $x\in A$ and $t\in (0,r/(4M))$,
\[
 |T_t(x)|\le |x|+ tM\le \frac{r}{2} +t M\le \frac{3r}{4}<r.
\]
Thus $T_t(A)\subset B_r$ i.e. $A\subset T_t^{-1}(B_r)$ as required.
\end{proof}
We end this section by recalling the formula for computing the density of the push-forward of absolutely continuous measures. We follow the presentation from \cite[Section 5.5]{ambrosio2005gradient}. To this aim we need to recall some definitions. For $f: \R^n\to \R^n$ we say that $z\in \R^n$ is the approximate limit of $f$ at $x$ if  the sets
\[
 \{y \ : \ |f(y)-z|>\eps\},
\]
have density $0$ at the point $x$. We then write $z=\tilde{f}$ for the approximate limit. Similarly, a linear map $L: \R^n\to \R^n$ is the approximate differential of $f$ at $x$ if  the sets
\[
 \lt\{y \ : \ \frac{|f(y)-\tilde{f}(x)-L(y-x)|}{|x-y|}>\eps\rt\},
\]
have density $0$ at the point $x$. We then write $L=\tilde{\nabla }f$. We let $\Sigma_f$ be the sets of points where $f$ is approximately differentiable. We then have, see \cite[Lemma 5.5.3]{ambrosio2005gradient}:
\begin{lemma}\label{lem:area}
 Assume that  $\rho\in L^1(\R^n)$ and that there exists $\Sigma\subset \Sigma_f$ such that $\tilde{f}\restr \Sigma$ is injective  and $|\{\rho>0\}\backslash \Sigma|=0 $. Then $f\# \rho
$ is absolutely continuous  if and only if $|\det \tilde{\nabla} f|>0$  a.e. on $\Sigma$. In this case, we have for a.e.  $x\in \Sigma$,
\[
 (f\#\rho)(\tilde{f}(x))= \frac{\rho(x)}{ |\det \tilde{\nabla} f(x)|}.
\]
\end{lemma}

\section{Two-point stability}\label{sec:2point}
The aim of this section is to prove Theorem \ref{twopointcompact}. We recall the definition \eqref{defLambda} of $\Lambda_\eps(\xi,f)$ and  will often write $\Lambda_\eps$ for $\Lambda_\eps(\xi,f)$ when there is no ambiguity.\\

Let us introduce some notation: for $r>0$, $\eta\ll1$, we set $\reps=r\eta$. For $\bar{x}\in (0,\reps)^n$  and $y\in \R^n$, we let $x_y$ be the projection of $y$ onto the grid $\reps\Z^n+ \bar z$ (which, for $\bar z$ fixed  is well-defined Lebesgue a.e.).  Setting $\X=\spt \lambda$, we let $X_{\bar z}=\{x\in \reps\Z^n+ \bar z \ : \ x=x_y \textrm{ for some }y\in \X\}$. We make the following assumptions on $\lambda$:
\begin{hypothesis}\label{hyp:lambda} \textcolor{white}{bla}
\begin{itemize}
 \item[(i)]  {there exists  $\eta\ll 1$ such that if $r\ll1$ there exists $M_0(r)\geq 1$ such that for every $\bar{z}\in (0,\reps)^n$ and $x, x' \in X_{\bar z}$} with  $|x-x'|\le \reps$ then
 \[\lambda(B_r(x)\cap B_r(x'))\ge  M_0(r)^{-1}\max(\lambda(B_r(x)),\lambda(B_r(x'))),\]
 \item[(ii)] $\lambda$ is absolutely continuous with respect to the Lebesgue measure,
 \item[(iii)] $X_{\bar z}$ is a connected subgraph of $\reps\Z^n$.
\end{itemize}
\end{hypothesis}
Notice that hypothesis (i) is a sort of doubling property for $\lambda$. It is for instance satisfied {with uniform $M_0$ }if $\lambda=\chi_\X$ for a set $\X$ having density estimates of the form $|\X\cap B_r(x)|\ges r^n$ for $x\in \X$, {as well as if $\lambda$ has a support with this property and a density bounded from above and below}. We assume that $\eps\ll1$ is such that (i) holds for $r=\eps$. All the implicit constants below depend on $\eta$.\\
 For each $\bar z\in (0,\reps)^n$ and  $x, x'\in X_{\bar z}$, we  let $\ell(x,x')$ be the (Euclidean) geodesic connecting $x$ to $x'$ in $\X+B_r$. We then let $\hat{\ell}(x,x')$ be an arbitrary curve in $X_{\bar z}\cap (\ell(x,x')+B_r)$ connecting $x$ to $x'$ with
 \[
  \HH^1(\hat{\ell}(x,x'))\les \HH^1(\ell(x,x')).
 \]
We let $I(x,x')=\hat{\ell}(x,x')\cap X_{\bar z}$ so that
\begin{equation}\label{lenghtI}
 r\HH^0(I(x,x'))\les  \HH^1(\hat{\ell}(x,x'))\les \HH^1(\ell(x,x')).
\end{equation}
For $p>1$, we then define the quantity (recall that $1/p+1/p'=1$)
\begin{multline}\label{deftaup}
 \tau_{p,\lambda}(r)=
 \\ \sup_{\bar z\in (0,\reps)^n}\sup_{z\in X_{\bar z}}\lt( \sum_{x,x'\in X_{\bar z}} \chi_{I(x,x')}(z)\lambda(B_r(x))\lambda(B_r(x')) \lt(\sum_{z'\in I(x,x')} \lambda(B_r(z'))^{-\frac{p'}{p}} \rt)^{\frac{p}{p'}}\rt).
\end{multline}
For $p = 1$, we instead set
\begin{multline*}\label{deftau1}
 \tau_{1,\lambda}(r)=
  \sup_{\bar z\in (0,\reps)^n}\sup_{z\in X_{\bar z}}\lt( \sum_{x,x'\in X_{\bar z}} \chi_{I(x,x')}(z)\lambda(B_r(x))\lambda(B_r(x')) \lt(\sup_{z'\in I(x,x')} \lambda(B_r(z'))^{-1} \rt)\rt).
\end{multline*}
When it is clear from the context we write $\tau(r)=\tau_{p,\lambda}(r)$. When $\bar z=0$ we simply write $X$ for $X_0$. Let us estimate $\tau_{p,\lambda}(r)$ under the assumptions of Theorem \ref{theo:introsimplequadratic}.
\begin{lemma}\label{lem:tauplambda}
 Let $\X$ be a Lipschitz domain and $\lambda\in \cP(\X)$ be such that $m\le \lambda\le M$ for some $m,M>0$. Then, for every $p\ge 1$ and $r\ll1$,
 \begin{equation}\label{eq:tauplambdaLip}
  \tau_{p,\lambda}(r)\les_\lambda r^{-(n+p-1)} \qquad \textrm{and }\qquad M_0(r)\les_\lambda 1.
 \end{equation}
The implicit constants depends on $\lambda$ only through $m, M$, $\diam(\X)$ and the Lipschitz constant of $\X$.
If instead $\X$ is convex then
\begin{equation}\label{eq:tauplambdaConv}
 \tau_{p,\lambda}(r)\les_\lambda r^{-p}.
\end{equation}

\end{lemma}
\begin{proof}
The estimate $M_0(r)\les_\lambda 1$ has already been discussed above so we focus only on the estimate of $\tau_{p,\lambda}$. We will consider only the case $p>1$; the case $p=1$ being analogous. Without loss of generality we may assume that $\bar z=0$. On the one hand, since  $\lambda(B_r(x))\approx_\lambda r^n$ for every $x\in X$,
\begin{multline*}
  \sum_{x,x'\in X} \chi_{I(x,x')}(z)\lambda(B_r(x))\lambda(B_r(x')) \lt(\sum_{z'\in I(x,x')} \lambda(B_r(z'))^{-\frac{p'}{p}} \rt)^{\frac{p}{p'}}\\
 \les_\lambda r^{n} \sum_{x,x'\in X} \chi_{I(x,x')}(z) \lt( \HH^0( I(x,x'))\rt)^{\frac{p}{p'}}.
\end{multline*}
Since $\X$ is Lipschitz, we have on the one hand $\HH^1(\ell(x,x'))\les_{\X} \diam(\X)$ for every $x,x'\in X$  and thus
\[
 \sup_{x,x'\in X} \HH^0(I(x,x'))\stackrel{\eqref{lenghtI}}{\les}r^{-1} \sup_{x,x'\in X} \HH^1(\ell(x,x')) \les_\X \diam(\X) r^{-1}.
\]
On the other hand, for every $z\in X$,
\begin{equation}\label{estimsumxx'Lip}
 \sum_{x,x'\in X}\chi_{I(x,x')}(z)\le (\HH^0(X))^2 \les_\X r^{-2n} |\X|^2.
\end{equation}
Combining these estimates yield \eqref{eq:tauplambdaLip}. Let us finally prove that if $\X$ is convex then \eqref{estimsumxx'Lip} may be improved to
\begin{equation}\label{estimsumxx'conv}
 \sum_{x,x'\in X}\chi_{I(x,x')}(z)\les_\X \diam(\X)^{n+1} r^{-(n+1)}.
\end{equation}
Set $R=\diam(\X)$ and notice that in this case $\ell(x,x')=[x,x']$ is a segment. Let $C_r(x,x')=\ell(x,x')+B_r$. We have by definition of $I$ and $\ell$,
\begin{multline*}
 \sum_{x,x'\in X}\chi_{I(x,x')}(z)\le \sum_{x,x'\in X}\chi_{C_r(x,x')}(z)\les_\X r^{-2n} \int_{\X\times\X} \chi_{C_r(x,x')}(z) dx dx'\\
 \le r^{-2n} \int_{\X} \int_{B_R(x)}\chi_{C_r(x,x')}(z) dx dx'.
\end{multline*}
Writing $x'=x+ t u$ for $u\in \partial B_1$ we find
\[
 \int_{\X} \int_{B_R(x)}\chi_{C_r(x,x')}(z) dx' dx\les R^{n-1}\int_{\X}\int_{\partial B_1}\int_0^R \chi_{C_r(x,x+ tu)}(z) dx d\HH^{n-1}(u) dt
\]
Finally notice that if $x\notin \R u +B_{2r}(z)$ then $z\notin  C_r(x,x+ tu)$ so that exchanging the order of integration we get
\begin{multline*}
 \int_{\X}\int_{\partial B_1}\int_0^R \chi_{C_r(x,x+ tu)}(z) dx d\HH^{n-1}(u) dt\le R\int_{\partial B_1} |\X\cap (\R u +B_{2r}(z))| d\HH^{n-1}(u) dt\\
 \les R^2 r^{n-1}.
\end{multline*}
This concludes the proof of \eqref{estimsumxx'conv}, from which we then obtain \eqref{eq:tauplambdaConv}.

% For the second estimate in \eqref{claimI} we observe that for every $z\in X$ we have
% \[
%  \sum_{x,x'\in X}\chi_{I(x,x')}(z)\le \lt(\sup_{x\in X} \sum_{x'\in X} \chi_{I(x,x')}(z) \rt) \#X.
% \]
% Since $\#X\les_\X r^{-n}|\X|$, it is enough to prove that for every $x,z\in X$,
% \[
%  \sum_{x'\in X} \chi_{I(x,x')}(z) \les_\X r^{-1} \diam(\X).
% \]
\end{proof}
\begin{remark}
 If we only assume that $\X$ is Lipschitz but not convex then \eqref{estimsumxx'Lip} is in general optimal, see Figure \ref{Figurestar}.
 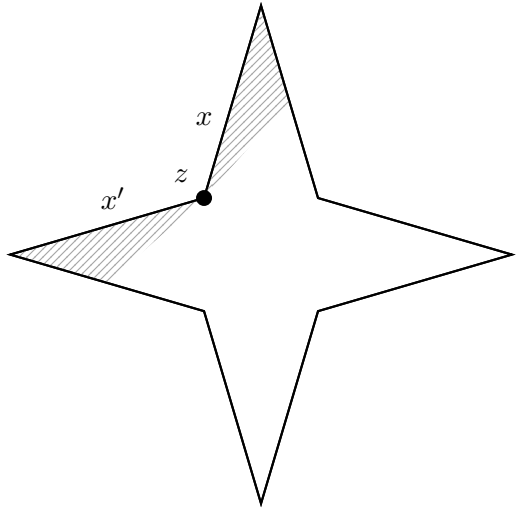
\begin{figure}[h]
\begin{tikzpicture}[scale=1.5, thick]
  % Coordinates for a sharper 4-pointed (ninja-style) star
  \coordinate (O) at (0,0);           % center
  \coordinate (A) at (0,2.2);         % top point
  \coordinate (B) at (2.2,0);         % right point
  \coordinate (C) at (0,-2.2);        % bottom point
  \coordinate (D) at (-2.2,0);        % left point

  % Inner corners pulled closer to center for sharper arms
  \coordinate (E) at (0.5,0.5);
  \coordinate (F) at (0.5,-0.5);
  \coordinate (G) at (-0.5,-0.5);
  \coordinate (H) at (-0.5,0.5);

  % Draw star outline
  \draw (A) -- (E) -- (B) -- (F) -- (C) -- (G) -- (D) -- (H) -- cycle;

  % Z point
  \fill (H) circle (2pt);
  \node at (-0.7,0.7) {$z$};
    % boundary of gray
     \coordinate (I) at (0.25,2.7/2);
     \coordinate (J) at (-2.7/2,-0.25);
  % Shaded triangular regions (top and left)
  \fill[pattern=north east lines, pattern color=gray!70] (H) -- (I) -- (A) -- cycle;
  \fill[pattern=north east lines, pattern color=gray!70] (H) -- (J) -- (D) -- cycle;

  % "x" labels inside shaded triangles
  \node at (-0.5,1.2) {$x$};
  \node at (-1.3,0.5) {$x'$};

  % Redraw outline for crisp edges
  \draw (A) -- (E) -- (B) -- (F) -- (C) -- (G) -- (D) -- (H) -- cycle;
\end{tikzpicture}
\caption{\label{Figurestar} For every point $x$ in the upper shaded region and every $x'$ in the bottom shaded region $z\in I(x,x')$.}
\end{figure}
\end{remark}

\begin{remark}\label{rem:taubounded}
 Notice that $\tau_{p,\lambda}(r)<\infty$ morally implies that $\X$ is bounded.  Assume for instance that   $\X$ contains the horizontal line. Since for every $x,x'$, and $p>1$,
\[
 \lt(\sum_{z'\in I(x,x')} \lambda(B_r(z'))^{-\frac{p'}{p}} \rt)^{\frac{p}{p'}}\ges \frac{1}{\lambda(B_r(x))}+\frac{1}{\lambda(B_r(x'))},
\]
we have for $z$ on the horizontal line,
\begin{multline*}
 \sum_{x,x'\in X} \chi_{I(x,x')}(z)\lambda(B_r(x))\lambda(B_r(x')) \lt(\sum_{z'\in I(x,x')} \lambda(B_r(z'))^{-\frac{p'}{p}} \rt)^{\frac{p}{p'}}\\
 \ges \sum_{x,x'\in X} \chi_{I(x,x')}(z)(\lambda(B_r(x))+\lambda(B_r(x')))=\infty.
\end{multline*}
\end{remark}

\begin{proof}[Proof of Theorem \ref{twopointcompact}]
 We only consider the case $p>1$ since the case $p=1$ is analogous. For every $x\in \R^n$, let
\[
 \m(x)=\int_{\R^n} |\xi(y)-f(x)|^p \rho_\eps(x-y) d\lambda(y),
\]
so that
\begin{equation*}\label{integlambda}
 \int_{\R^n} \m(x)dx=\Lambda_\eps.
\end{equation*}
Let $\eta\ll1$  be such that point (i) of Hypothesis \ref{hyp:lambda} holds for $r\le \eps$, and set {$\reps=\eta r$,  $Q_{\reps}=(0, \reps)^n$ and for $z\in \reps \Z^n$, $Q_{\reps}(z)=Q_{\reps}+z$.}  We then have
\[
 \int_{Q_{\reps}} \sum_{x\in  \reps\Z^n} \m(x+z)dz=\sum_{x\in  \reps\Z^n}\int_{Q_{\reps}(x)} \m(z) dz\le \int_{\R^n} \m(z) dz =\Lambda_\eps,
\]
so there exists $\bar{z}\in Q_{\reps}$ such that
\begin{equation}\label{estimsumlambda}
 {r^{n}\sum_{x\in  \reps\Z^n} \m(x+\bar{z})\les  \Lambda_\eps.}
\end{equation}
Up to a translation, we can assume that $\bar{z}=0$. We estimate
\begin{multline*}
 \int_{\R^n\times \R^n}|\xi(y)-\xi(y')|^p d\lambda(y) d\lambda(y')\les \int_{\R^n\times \R^n}|\xi(y)-f(x_y)|^p d\lambda(y) d\lambda(y')\\+\int_{\R^n\times \R^n}|f(x_y)-f(x_{y'})|^p d\lambda(y) d\lambda(y')
 +\int_{\R^n\times \R^n}|\xi(y')-f(x_{y'})|^p  d\lambda(y) d\lambda(y')\\
 \les \int_{\R^n}|\xi(y)-f(x_{y})|^p d\lambda(y) +\int_{\R^n\times \R^n}|f(x_{y})-f(x_{y'})|^p  d\lambda(y) d\lambda(y').
\end{multline*}
We start by considering the first right-hand side term. Since $|y-x_y|\le \sqrt{n}/2 \reps<\eps$,  we have $y\in B_{\eps}(x_y)$ and thus as $1 \les \eps^{n} \rho_\eps$ on $B_\eps$ (recall that we assumed $\rho \ges 1$ on $B_1$), so
\begin{multline*}
 \int_{\R^n}|\xi(y)-f(x_{y})|^p d\lambda(y)= \sum_{x\in \reps\Z^n} \int_{\{y \ :\ x_y=x\}} |\xi(y)-f(x)|^p d\lambda(y)\\
 \le \sum_{x\in \reps\Z^n} \int_{B_\eps(x)} |\xi(y)-f(x)|^p d\lambda(y)
 \les \eps^n \sum_{x\in  \reps\Z^n} \m(x)\stackrel{\eqref{estimsumlambda}}{\les} {(\eps/r)^n\Lambda_\eps.}
\end{multline*}
 For the second sum, we can similarly write
\[
 \int_{\R^n\times \R^n}|f(x_{y})-f(x_{y'})|^p  d\lambda(y) d\lambda(y')\les \sum_{x,x'\in X} \lambda(B_r(x))\lambda(B_r(x')) |f(x)-f(x')|^p.
\]
To lighten notation a bit, set $\lambda_x=\lambda(B_r(x))$. Let now $x,x'\in X$. If $x$ and $x'$ are nearest neighbors in $\reps\Z^n$, by Hypothesis (i), letting $A= B_r(x)\cap B_r(x')$ we have $\lambda(A)\ge M_0(r)^{-1}\max(\lambda_x,\lambda_{x'})$. By the triangle inequality and Jensen we then find
\begin{multline*}
 |f(x)-f(x')|^p\les \lt|f(x) -\frac{1}{\lambda(A)}\int_{A} \xi d\lambda\rt|^p+ \lt|f(x')-\frac{1}{\lambda(A)}\int_A \xi d\lambda\rt|^p \\
 \le \frac{1}{\lambda(A)}\int_{A} |f(x)-\xi(y)|^p d\lambda(y)+  \frac{1}{\lambda(A)}\int_{A} |f(x')-\xi(y)|^p d\lambda(y)\\
 \le \frac{1}{\lambda(A)}\int_{B_r(x)} |f(x)-\xi(y)|^p d\lambda(y)+  \frac{1}{\lambda(A)}\int_{B_r(x')} |f(x')-\xi(y)|^p d\lambda(y)\\
 \les \eps^n M_0(r)\lt(\frac{\m(x)}{\lambda_x}+\frac{\m(x')}{\lambda_{x'}}\rt).
\end{multline*}
Take now $x, x'\in X$ arbitrary. Then by the triangle inequality and the above computation, summing along the geodesic we get
\begin{multline*}
 |f(x)-f(x')|^p \les \eps^n M_0(r) \lt( \sum_{z\in I(x,x')} \lt(\frac{\m(z)}{\lambda_z}\rt)^{\frac{1}{p}}\rt)^p\\
 \le \eps^n M_0(r)\lt(\sum_{z\in I(x,x')} \lambda_z^{-\frac{p'}{p}} \rt)^{\frac{p}{p'}} \sum_{z\in I(x,x')} \m(z),
\end{multline*}
where the last inequality is due to the H\"older inequality. After summation, we find
\begin{multline*}
  \sum_{x,x'\in X} \lambda_x\lambda_{x'} |f(x)-f(x')|^p
  \les \eps^n M_0(r) \sum_{x,x'\in X} \lambda_x\lambda_{x'} \lt(\sum_{z\in I(x,x')} \lambda_z^{-\frac{p'}{p}} \rt)^{\frac{p}{p'}} \sum_{z\in I(x,x')} \m(z) \\
 =\eps^n M_0(r)\sum_z \lt( \sum_{x,x'\in X} \chi_{I(x,x')}(z)\lambda_x\lambda_{x'} \lt(\sum_{z'\in I(x,x')} \lambda_{z'}^{-\frac{p'}{p}}\rt)^{\frac{p}{p'}}\rt) \m(z).
\end{multline*}
By definition \eqref{deftaup} of $\tau_{p,\lambda}$, we get
\[
 \sum_{x,x'\in X} \lambda_x\lambda_{x'} |f(x)-f(x')|^p
  \les M_0(r)\tau_{p,\lambda}(r) \eps^n \sum_z \m(z)\stackrel{\eqref{estimsumlambda}}{\les} {M_0(r)(\eps/r)^n\tau_{p,\lambda}(r) \Lambda_\eps.}
\]
In conclusion, we find
\begin{multline*}\label{estimdoubleintegxi}
 \int_{\R^n\times \R^n}|\xi(y)-\xi(y')|^p d\lambda(y) d\lambda(y')\les (\eps/r)^n\Lambda_\eps + M_0(r)(\eps/r)^n\tau_{p,\lambda}(r) \Lambda_\eps\\
 \les M_0(r)(\eps/r)^n\tau_{p,\lambda}(r) \Lambda_\eps.
\end{multline*}
Recalling that
\[z=\int_{\R^n} \xi d\lambda, \]
we find by Jensen's inequality,
\[
 \int_{\R^n}|\xi-z|^p d\lambda\leq \int_{\R^n\times \R^n}|\xi(y)-\xi(y')|^p d\lambda(y) d\lambda(y')\les {M_0(r)(\eps/r)^n\tau_{p,\lambda}(r) \Lambda_\eps,}
\]
which concludes the proof of \eqref{eq:twopointcompact}.\\
\end{proof}
\begin{remark}\label{rmk:BBM}
    Note that if we assume that $\X$ is convex, as a consequence of Theorem \ref{twopointcompact} and Lemma \ref{lem:tauplambda} (and in particular of \eqref{eq:tauplambdaConv}), we can obtain a new characterization of  constant functions: If $\xi$ is  such that, for every $\eps > 0$ there exists a function $f_\eps$ such that $\Lambda_\eps(\xi, f_\eps) = o(\eps ^ p)$, then $\xi$ is constant. This may be seen as a generalization of \cite{brezis2002recognize}.
    
    For another consequence, we can take $\lambda$ to be the uniform probability measure on $\X$ and $\xi = f$. {To lighten notation, let us take $\eps = r$ and omit the dependence on $M_0$.} We get
    \[
    \int_\X |f(x) - z| ^ p dx \les \tau_{p,\lambda}(\eps)\Lambda_\eps(f,f) \les_{\X} \eps ^ {-p} \int_{\X \times \R ^ n} \rho_\eps(x-y)|f(x)-f(y)|^p dx dy.
    \]
    Since $|z|\les \eps$ on $\spt \rho_\eps$,  we obtain the estimate:
    \[
    \int_{\X} |f-z|^p dx\les_{\X} \int_{\X \times \R ^ n} \rho_\eps(x-y) \frac{|f(x)-f(y)|^p}{|x-y|^p} dx dy.
    \]
    This estimate is similar to the one proved in \cite{Ponce04}. It is a non-local version  of the Poincar\'e-Wirtinger inequality which actually implies the standard one  as a consequence of \cite[Theorem 1]{Bourgain01}.
\end{remark}
We now show that a similar estimate to that of Theorem \ref{twopointcompact} can be obtained for some non-compactly-supported $\lambda$. For simplicity we restrict ourselves to the Gaussian case, see also Remark \ref{rem:logconcavProp34}.
\begin{proposition}\label{prop:twopointgaussian}
 Let $\lambda$ be a standard Gaussian of mean zero and variance $1$. Assume also that there exists a constant $C>0$ such that
 \begin{equation}\label{hypxigauss}
  |\xi(x)|\le C( 1+|x|).
 \end{equation}
Then, with
\[
 z=\int_{\R^n} \xi d\lambda,
\]
we have for every $f$, $\eps \in (0,1)$, $p\ge 1$, and any $\beta>0$
\begin{equation}\label{eq:twopointgaussian}
\int_{\R^n} |\xi-z|^pd\lambda \les_{\rho, C, \beta} {\eps^{-p}} \Lambda_\eps(\xi,f)^{1-\beta}.
\end{equation}
\end{proposition}
\begin{proof}
We consider only the case $p>1$ since the case $p=1$ is analogous. Notice first that by \eqref{hypxigauss},
\[
 \int_{\R^n} |\xi-z|^p d\lambda\les_C 1,
\]
and thus we may assume without loss of generality that $\Lambda_\eps\ll1$.
 Let  $z=\int_{\R^n} \xi d\lambda$ and for $R\gg1$ to be chosen,
\[
 z_R=\int_{B_R} \xi d\lambda, \qquad \textrm{and } \qquad \lambda_R= \lambda \chi_{B_R}.
\]
We then write
\begin{equation}\label{decompxiz}
 \int_{\R^n} |\xi-z|^pd\lambda\les \int_{\R^n} |\xi -z_R|^p d\lambda_R + |z-z_R|^p + \int_{B_R^c} |\xi-z|^p d\lambda.
\end{equation}
Let us first estimate the second right-hand side term in \eqref{decompxiz}. Since $\lambda$ is radially symmetric, we have
{$$\int_{B_R^c} |x| d\lambda = C(n)\int_R^\infty r^ne^{-r^2/2}dr \les R^{n-1}e^{-R^2/2}$$
and therefore}
\begin{equation}\label{zzR}
 |z-z_R|^p=\lt|\int_{B_R^c} \xi d\lambda\rt|^p\stackrel{\eqref{hypxigauss}}{\les_C}\lt(\int_{B_R^c} |x| d\lambda\rt)^p\les (R^{n-1} \lambda(R))^p.
\end{equation}
We now turn to the third right-hand side term in \eqref{decompxiz}. By \eqref{hypxigauss} we have similarly
\begin{equation}\label{xiRlarge}
 \int_{B_R^c} |\xi-z|^p d\lambda\les \int_{B_R^c} |\xi|^p d\lambda+ |z|^p \lambda(B_R^c)\les_C \int_{B_R^c} |x|^p d\lambda\les R^{n-2+p} \lambda(R).
\end{equation}
Let us point out that since $p\ge 1$, for $R$ large the right-hand side of \eqref{xiRlarge} is larger than the right-hand side of \eqref{zzR}. We finally estimate the first right-hand side term in \eqref{decompxiz}. By Theorem \ref{twopointcompact} applied to $\lambda_R/\lambda(B_R)$ and setting $\tau_R(r)=\tau_{p,\lambda_R}(r)$, we have for $r\le \eps$
\begin{equation}\label{firstgauss}
 \int_{\R^n} |\xi-z_R|^p d\lambda_R\les \tau_R(\eps) \int |\xi(y)-f(x)|^p\rho_\eps(x-y) dx d\lambda_R(y)\les {M_0(r)(\eps/r)^n\tau_R(r)\Lambda_\eps},
\end{equation}
so we are left with estimating $\tau_R(r)$ {and $M_0(r)$} for an appropriate choice of $r$. We first note that { since $\lambda_R/\lambda(B_R)$ has nice support, part $(i)$ of Hypothesis \ref{hyp:lambda} holds with
\begin{align*}\label{M}
 M_0(r)&\les \sup_{|x-x'| \leq r}\frac{\sup_{u \in B_R \cap(B_r(x) \cap B_r(x'))} \lambda_R(u)}{\inf_{u \in B_R \cap(B_r(x) \cap B_r(x'))} \lambda_R(u)} \notag \\
&\les \exp(4(rR + r^2)).
\end{align*}
Hence if we take $r = \min(\eps, R^{-1})$ we have $M_0(r) \les 1$. }

{We set
\begin{equation}\label{choiceR}
 R=\sqrt{2\log \Lambda_\eps }.
\end{equation}
Since $r \le R^{-1}$, we have for every $x\in X$,
\[
 \lambda(B_r(x))\approx \lambda(x) |B_r|,
\]}
and thus for every $(x,x')\in X$, assuming $|x'|>|x|$:
\begin{multline*}
 \lt(\sum_{z'\in I(x,x')} \lambda(B_r(z'))^{-\frac{p'}{p}}\rt)^{\frac{p}{p'}} \approx |B_r|^{-1} \lt(\sum_{z'\in I(x,x')} \frac{1}{\lambda(z')^{\frac{p'}{p}}}\rt)^{\frac{p}{p'}}\\
 \les |B_r|^{-1} \lt( r^{-1} \int_{-|x|'}^{|x'|} \exp\lt(\frac{p'}{2p} t^2\rt)\rt)^{\frac{p}{p'}}\\
 \les |B_r|^{-1} \lt( r^{-1} \min(1,|x'|^{-1})\rt)^{\frac{p}{p'}}\lambda^{-1}(x').
 %\simeq r^{-(d-1)} \frac{|x|}{\lambda(x)}+ r^{-(d-1)} \frac{|x'|}{\lambda(x')}.
\end{multline*}
We thus find
\begin{multline*}
 \sum_{x,x'\in X} \chi_{I(x,x')}(z) \lambda(B_r(x))\lambda(B_r(x'))\lt(\sum_{z'\in I(x,x')} \lambda(B_r(z'))^{-\frac{p'}{p}}\rt)^{\frac{p}{p'}} \\
 \les r^{n-\frac{p}{p'}} \sum_{x,x'\in X} \chi_{I(x,x')}(z) \lambda(x) \min(1,|x'|^{-\frac{p}{p'}})\\
 \le r^{n-\frac{p}{p'}} \sum_{x,x'\in X} \chi_{I(x,x')}(z) \lambda(x).
\end{multline*}
Using the same  computations as in the proof of \eqref{estimsumxx'conv} from Lemma \ref{lem:tauplambda} we have
\begin{multline*}
 \sum_{x,x'\in X} \chi_{I(x,x')}(z) \lambda(x)\les r^{-2n} R^{n} \int_{\partial B_1} \lambda(B_R\cap (\R u+B_{2r}(z))) d\HH^{n-1}(u)
 \les r^{-(n+1)} R^n.
\end{multline*}
From this we get
\[
 \tau_R(r)\les r^{-p} R^n.
\]
Plugging this into \eqref{firstgauss} and recalling that $M_0(r) \les 1$ yields
\begin{equation}\label{firstgaussfinal2}
  \int_{\R^n} |\xi-z_R|^p d\lambda_R\les \eps^n r^{-(p+n)} R^n   \Lambda_\eps.
\end{equation}
Plugging \eqref{zzR}, \eqref{xiRlarge} and \eqref{firstgaussfinal2} into \eqref{decompxiz} leads to
\[
 \int_{\R^n} |\xi-z|^pd\lambda\les \eps^n r^{-(p+n)} R^n   \Lambda_\eps+R^{n-2+p}  \exp(-R^2/2).
\]
By the choice  \eqref{choiceR} of $R$ and since $r = \min(\eps, R^{-1}) \geq \eps R^{-1}$, we find
\[
 \int_{\R^n} |\xi-z|^pd\lambda\les \lt(\eps^{-p} |\log \Lambda_\eps|^{\frac{2n+p}{2}} + |\log \Lambda_\eps|^{n-2+p}\rt) \Lambda_\eps,
\]
which concludes the proof of \eqref{eq:twopointgaussian}.
\end{proof}
%\begin{remark}
% Notice that in order to have $\eps R\ll 1$, we need to have
%\[
% |\log \Lambda_\eps|\gg \eps^{-2} \qquad \Longleftrightarrow \qquad \Lambda_\eps\les \exp(-C \eps^{-2}),
%\]
%which is extremely small.
%\end{remark}
\begin{remark}\label{rem:logconcavProp34}
 Since we used bounds from above and below on the measure of balls, it is not clear if Proposition \ref{prop:twopointgaussian} holds for more general uniformly log-concave measures.
\end{remark}

\section{Stability of Potentials for Linear Cost}\label{sec:stablin}
The aim of this section is to prove Theorem \ref{thm1}. We will also obtain some extensions for measures which  have densities not bounded from below, see Theorem \ref{thm2}. This includes also some measures with unbounded support.

For $\lambda$ absolutely continuous with respect to the Lebesgue measure with $\X=\spt \lambda$ and $\mu$ a measure of same mass, we let $T$ be an optimal transport map for $W_1(\lambda,\mu)$ and $\psi$ be an optimal Kantorovich potential. Recall the definitions \eqref{defTt} of $T_t$, \eqref{defmuj} of $\mu_t$, $j_t$ and \eqref{eq:transportdens} of $\sigma$. In order to prove Theorem \ref{thm1} we first obtain a stability inequality when integrating with respect to $\sigma$.
\begin{lemma}\label{lem:centralstab1}
     For every $1$-Lipschitz function $\phi$,
    \begin{equation}\label{eq:stabsigma}
        \|\nabla \psi - \nabla \phi\|_{L ^ 2(\sigma)} ^ 2 \leq 2 \int_{\R^n} (\psi - \phi) d(\mu - \lambda).
    \end{equation}
\end{lemma}
\begin{proof}
    Using $\mu=T\#\lambda$ and $\psi(T(x))-\psi(x)=|T(x)-x|$, we calculate:
    \[
        \begin{aligned}
            \int_{\R^n} (\psi - \phi) d(\mu - \lambda) & = \int_{\R^n} (|T(x)-x| - \phi(T(x)) + \phi(x)) d\lambda \\
            & = \int_{\R^n} \int_0 ^ 1 (|T(x)-x| - \ang{\nabla \phi(T_t(x)), T(x)-x}) dtd\lambda \\
            & = \int_{\R^n} \int_0 ^ 1 |T(x)-x| \left(1 - \ang{\nabla \phi(T_t(x)), \frac{T(x)-x}{|T(x)-x|}}\right) dtd\lambda.
            \end{aligned}
            \]
Since $|\nabla \phi|\le 1$ we have
\[
 1 - \ang{\nabla \phi(T_t(x)), \frac{T(x)-x}{|T(x)-x|}}\ge \frac{1}{2}\lt|\nabla \phi(T_t(x))-\frac{T(x)-x}{|T(x)-x|}\rt|^2.
\]
Using moreover that
\[
 \frac{T(x)-x}{|T(x)-x|}=\nabla \psi(T(x))=\nabla \psi(T_t(x)),
\]
we find
            \[\begin{aligned}
             \int_{\R^n} (\psi - \phi) d(\mu - \lambda)& \geq  \frac{1}{2} \int_0 ^ 1 \int_{\R^n} |T(x)-x| |\nabla \phi(T_t(x)) - \nabla \psi(T_t(x))| ^ 2 d\lambda dt \\
            &\stackrel{\eqref{defmuj}}{ =} \frac{1}{2} \int_0 ^ 1\int_{\R^n} |\nabla \psi - \nabla \phi| ^ 2 d|j_t|dt \\
            & \stackrel{\eqref{eq:transportdens}}{=} \frac{1}{2} \|\nabla \psi - \nabla \phi\|_{L ^ 2(\sigma)} ^ 2.
        \end{aligned}
    \]
\end{proof}
In order to prove Theorem \ref{thm1}, we thus need to bound $\sigma$ from below by $\lambda$. This may be seen as a quantitative version of Lemma \ref{lem:qualitabscont}.  Notice that bounds of  $\sigma$ from above by $\lambda$ are available in the literature, see \cite{de2002regularity} or \cite[Section 4.3]{Santam}. As we will see the bound will degenerate close to $\partial \X$. We will thus use the following elementary (and probably well-known) result. For $\X\subset \R^n$  compact and $r>0$ we define the $r-$erosion
\begin{equation}\label{Xr}
\X_r = \{x \in \X : d(x, \partial \X) \geq r\}.\end{equation}

\begin{lemma}\label{lem:Minkowski}
    Let $\X\subset \R^n$ be a compact set such that $\partial \X$ is  rectifiable.  Then for every $\alpha \in (0, 1)$,
    \begin{equation}\label{Ifin}
    \begin{aligned}
    I_\alpha(\X) = \int_\X d(x, \partial \X) ^ {-\alpha} dx<\infty.
    \end{aligned}
    \end{equation}
    Moreover, if  $\partial \X$ is $C^2$  then
    \begin{equation}\label{limsupI}
     \limsup_{r\to 0} I_\alpha(\X_r)<\infty.
    \end{equation}
Finally, if
    $\X$ is  a convex set with $ \X\subset B_R$ for some $R>0$ then
    \begin{equation}\label{MinkowskiConvex}
      I_\alpha(\X) \les_\alpha R^{n-1}.
    \end{equation}

\end{lemma}
\begin{proof}
     We start with the general case. By definition of upper Minkowski content \cite[Definition 2.100]{AFP} and \cite[Theorem 2.106]{AFP}, there exists $r_0>0$ such that for $r\le r_0$,
    \begin{equation}\label{difvolrectif}
     |\X \backslash \X_r|\les_\X r.
    \end{equation}
    Therefore,
    \[
    \begin{aligned}
        I_\alpha(\X) = \int_\X d(x, \partial \X) ^ {-\alpha} dx & = \int_0 ^ \infty |\{x \in \X: d(x, \partial \X) ^ {-\alpha} \geq s\}| ds \\
        & \leq \int_0 ^{r_0 ^ {-\alpha}} |\X| ds + \int_{r_0 ^ {-\alpha}}^\infty |\X \backslash \X_{(1 / s) ^ {1 / \alpha}}| ds \\
        & \les_{\X} 1 +  \int_{r_0 ^ {-\alpha}} ^ \infty s^{-1/\alpha}ds \\
        & <\infty.
    \end{aligned}
    \]
    For every compact set $\X$, denoting $f(x)=d(x,\partial \X)$ which satisfies $|\nabla f|=1$ a.e. and $\X^{(r)}=\{ x\in \X \ : d(x, \partial \X) = r\}$ we have by the co-area formula,
    \[
     |\X \backslash \X_r|=\int_{\X\backslash \X_r} |\nabla f| dx=\int_0^r \mathcal{H}^{n-1}(\X^{(s)})ds\le \lt(\sup_{s\in[0,r]} \mathcal{H}^{n-1}(\X^{(s)})\rt) r.
    \]
If now $\partial \X$ is $C^2$ then on the one hand, for $r_0$ small enough, see e.g. \cite{AmbDanc},
\[
 \lt(\sup_{s\in[0,r_0]} \mathcal{H}^{n-1}(\X^{(s)})\rt)\les_\X 1
\]
and on the other hand, for $\eta,s$ small enough $(\X_\eta)^{(s)}=\X^{(s+\eta)}$ so that for $\eta, r$ small enough we have a uniform in $\eta$ version of \eqref{difvolrectif},
\[
 |\X_\eta \backslash (\X_\eta)_r|\les_\X r.
\]
Injecting this in the previous proof yields \eqref{limsupI}.\\
    If finally $\X$ is convex, then $\X_r$ is also convex with $\X_r\subset \X\subset B_R$
so that $\mathcal{H}^{n-1}(\partial \X_r)\le \mathcal{H}^{n-1}(\partial B_R)\les R^{n-1}$ and since $\partial \X_r= \X^{(r)}$, this yield the quantitative version of \eqref{difvolrectif}
 \[
 |\X\backslash \X_r|\les R^{n-1} r.
 \]
    Plugging this back in the previous computation (with the choice $r_0=1$) concludes the proof of \eqref{MinkowskiConvex}.
\end{proof}
 We can use this lemma to obtain a relation between the measure $\lambda$ and the transport density $\sigma$, namely an upper bound on the Rényi divergence $D_\alpha(\lambda || \sigma)$. Recall that for arbitrary probability measures $\rho_1, \rho_2$ in $\mathbb R ^ d$, the Rényi divergence is given by
\[
D_\alpha(\rho_1 || \rho_2) = \frac{1}{\alpha - 1}\ln\left(\int_{\mathbb R ^ n} \left(\frac{d\rho_1}{d\rho_2}\right) ^ {\alpha - 1} d\rho_1\right).
\]
{See \cite{bobkov2025} and references therein for background on R\'enyi divergences and comparisons with various classical distances between probability measures.}
\begin{lemma}\label{lem:Renyi}
    Under the assumptions of Theorem \ref{thm1}, the following inequality holds:
    \be
    \label{eq:Renyi}
        D_\alpha(\lambda || \sigma) \leq C_{\X, \alpha, R} + \frac{\alpha}{\alpha - 1}\ln(M) - \ln(m),
    \ee
    for every $\alpha < 3 / 2$. In particular we have $\lambda\ll \sigma$.

    Moreover, the constant $C_{\X, \alpha, R}$ in \eqref{eq:Renyi} can be chosen to be of the form
    \[
        C_{\X, \alpha, R} = \frac{1}{\alpha - 1} \ln(c_1 |\X| + c_2I_{2(\alpha - 1)}(\X)),
    \]
    where $c_1, c_2$ depend only on $\alpha, R$ and the dimension $n$.
\end{lemma}
\begin{proof}
    We first prove the lemma under the additional assumptions that $\mu$ is discrete, and that $n \geq 2$. Let $T$ be an optimal transport map from $\lambda$ to $\mu$, and denote as usual, $T_t = tT + (1 - t)\mathrm{id}$. We recall the formula \eqref{eq:transportdens} for the transport density $\sigma$. From now on, we assume that $t \neq 0, 1$. As a byproduct of the proof in \cite[Theorem 4.16]{Santam}, the map $T_t$ is essentially injective in this case. Writing $\mu=\sum_{i=1}^I \mu_i \delta_{z_i}$ for some $\mu_i>0$ and $z_i\in \R^n$, we set for every $i=1,\cdots, I$
    \[
     A_i=\{x \ : \ T(x)=z_i\}.
    \]
    We claim that for every point $x$ of density one of $A_i$, the map $T_t$ is approximately differentiable at $x$ with (recall the notation from Lemma \ref{lem:area} for approximate limits and differentials)
    \begin{equation}\label{claimtildeTt}
    \tilde{T_t}(x)= (1-t)x+ tz_i\qquad \textrm{ and }\qquad \tilde{\nabla}T_t(x)=(1-t) \mathrm{id}.\end{equation}
    It is of course enough to prove the corresponding claim for $T$ instead of $T_t$. Let $x$ be a point of $A_i$ of density one. Since $\mu$ is discrete we have for $\eps>0$ small enough,
    \[
     \{ y \ : \ |T(y)-z_i|>\eps\}= A_i^c
    \]
which indeed has density zero at $x$. This proves that $\tilde{T}(x)=z_i$. Similarly, {since the support of $\lambda$ is bounded}
\[
 \{ y \ : |T(y)-z_i|/|y-x|>\eps\}\subset A_i^c
\]
which gives as claimed $\tilde{\nabla} T(x)=0$.\\
By \eqref{claimtildeTt} we have for a.e. $x\in \X$, $0<\det \tilde{\nabla}T_t(x)=(1-t)^n<1$ and thus by Lemma \ref{lem:area}, for a.e. $x\in \X$,
    \[
        |j_t|(T_t(x)) = |T(x) - x| \lambda(x) / |\det \tilde{\nabla} T_t(x)| > |T(x) - x| \lambda(x)\ge |T(x)-x|m,
    \]
    where we used that $\lambda\ge m$ on $\X$ by hypothesis. Let now $y\in T_t(\X)\cap \X$ and set $x= T_t^{-1}(y)$. Note that since $T_t$ is essentially injective, $T_t ^ {-1}$ is indeed well-defined for almost all such $y$. Since $y$ lies in the ray between $x$ and $T(x)$ we have $|T(x)-x|\ge |T(x)-y|$. Moreover, since $y\in \X$ and $T(x)\in \X^c$ by the hypothesis that $\spt \lambda\cap \spt \mu=\emptyset$, we have $|T(x)-y|\ge d(y,\partial \X)$ and  the  previous inequality implies the lower bound
    \be
    \label{lambdatlowerbound}
     |j_t|(y)\ge d(y,\partial \X) m \qquad \textrm{for } a.e.  \ y\in T_t(\X)\cap \X.
    \ee

    For $y \in \X$, let
    \[t(y)=\sup\{ t_0 \ : \ y\in T_t(\X)\cap \X \textrm{ for every } t<t_0\}.\]
    Then, according to \eqref{eq:transportdens} and \eqref{lambdatlowerbound}, we have for a.e. $y\in \X$,
    \be
    \label{sigmalowerbound}
        \sigma(y) = \int_0 ^ 1 |j_t|(y) dt \geq \int_0 ^ {t(y)} |j_t|(y) dt > \int_0 ^ {t(y)} d(y, \partial \X) m dt = t(y) d(y, \partial \X) m.
    \ee
    Notice that since $\sigma$ is absolutely continuous, see \cite[Theorem 4.16]{Santam}, we may use Fubini to exchange the integral in time and the a.e. statement \eqref{lambdatlowerbound}.  We now want to find a relatively large subset of $\X$ for which the right-hand side of \eqref{sigmalowerbound} can be bounded from below.
    Recall the definition  \eqref{Xr} of $\X_r$,  that the supports of $\lambda, \mu$ both lie in a ball of radius $R$ and  define for $t_0\in (0,1)$,
    \[
        \X(t_0) = T_{t_0}\left(\X_{t_0(2R + 1)}\right).
    \]
    We use the convention that $\X(t_0)=\emptyset$ if $\X_{t_0(2R + 1)}=\emptyset$.
    We claim that for almost every $y \in \X(t_0)$, we have $y\in \X_{t_0}$, and $t(y) \geq t_0$. By \eqref{sigmalowerbound} this would imply, that
    \be
    \label{sigmaImprovedLowerBound}
        \sigma(y) > t(y) d(y, \partial \X) m \geq mt_0 ^ 2 \qquad \qquad  \textrm{on } \ \X(t_0).
    \ee
Notice that these claims are empty if $\X_(t_0)$ is empty.
    For the first claim, let $x \in \X_{t_0(2R + 1)}$. The points $x, T(x)$ both lie in a ball of radius $R$, so $|T(x) - x| \leq 2R$. If we denote $y = T_{t_0}(x)$, we see that
    \[
        |y - x| = |T_{t_0}(x) - x| = t_0|T(x) - x| \leq 2t_0R,
    \]
    so $y \in \X_{t_0}$, since $ d(y, \partial \X)\ge d(x, \partial \X)- |x - y| { \geq t_0(1+2R) - 2Rt_0}$.
    
    The second claim will require more effort. Let $A$ be the set of all $x \in \X$ which lie on a line connecting two points in the support of $\mu$. Since $\mu$ is discrete and $n\geq 2$, $A$
    is equal to a finite union of lines, and thus has {Lebesgue} measure zero. Notice that if $x\in \X\backslash A$ and $L$ is a line containing $x$ then $\HH^0(L\cap \spt \mu)\le 1$. Let $y=T_{t_0}(x) \in \X(t_0) \backslash T_{t_0}(A)$. For $t < t_0$, define $z_t = \frac{y - tT(x)}{1 - t}$. We aim to show that $T_t(z_t) = y$, which would imply the second claim. { Note that $T(z_t)$ is well-defined since $d(x,\partial \X) \geq t_0(1+2R)$ and
		$$d(z_t, x) = \frac{t_0-t}{1-t}|x-T(x)| \leq 2R t_0$$
		so that $z_t \in \X$. }
    
    By definition $z_t$ lies between $x$ and $T(x)$. Recall that transport rays are one-dimensional segments covering our entire domain, and that $x, T(x)$ always lie on the same transport ray. It is well-known \cite[Corollary 3.8]{Santam} that transport rays only intersect at their endpoints, so the only transport ray that $z_t$ belongs to is the one containing $x, T(x)$. Since $z_t, T(z_t)$ also lie on the same transport ray, we conclude that $x, T(x), T(z_t)$ are all on the same line $L$. However since $x\in \X\backslash A$ and $T(z_t),T(x)\in L\cap \spt \mu$,  we must have $T(z_t) = T(x)$. Therefore,
    \[
        T_t(z_t) = tT(z_t) + (1 - t)z_t = tT(x) + (y - tT(x)) = y.
    \]
    This concludes the proof of \eqref{sigmaImprovedLowerBound}.
     We estimate now the size of $\X(t_0)$ in relation to $\X$. Since the map $T_t$ is essentially injective and $(1 - t)$-contractive, we see  that
    \[
        |\X(t_0)| = (1 - t_0) ^ n |\X_{t_0(2R + 1)}| \geq (1 - nt_0) |\X_{t_0(2R + 1)}|.
    \]
    Therefore, the size of the set $\X \backslash \X(t_0)$ can be bounded as
    \[
        \begin{aligned}
        |\X \backslash \X(t_0)| & \leq |\X| - (1 - nt_0)|\X_{t_0(2R + 1)}| = |\X| - |\X_{t_0(2R + 1)}| + nt_0|\X_{t_0(2R + 1)}| \\
        & \leq \left|\X \backslash \X_{t_0(2R + 1)}\right| + n|\X|t_0.
        \end{aligned}
    \]
    Notice that this estimate trivially holds if $\X_{t_0(2R + 1)}=\emptyset$.
    We now prove  that the integral of $\sigma ^ {-\beta}$ converges, as long as $\beta < 1 / 2$. To simplify notation, we denote $\sigma' = \sigma / m$. By \eqref{sigmaImprovedLowerBound}, for almost all $y$ and $s>0$, if $\sigma'(y) \leq s$ then $y\not\in \X(s ^ {1 / 2})$. Therefore,
    \[
        \begin{aligned}
        \int_\X (\sigma') ^ {-\beta} dx & = \int_0 ^ \infty \left|\left\{y : \sigma'(y) < s ^ {-1 / \beta}\right\}\right| ds \\
        & \leq \int_0 ^ \infty \left|\X \backslash \X\left(s ^ {-1 / 2\beta}\right)\right| ds \\
        & \leq \int_0 ^ 1 |\X| ds + \int_1 ^ \infty \left|\X \backslash \X_{s ^ {-1 / 2\beta}(2R + 1)}\right| ds + \int_1 ^ \infty n|\X| s ^ {-1 / 2\beta} ds \\
        & = \left(\frac{2\beta}{1-2\beta} n + 1\right)|\X| + \int_1 ^ \infty \left|\left\{x : d(x, \partial \X) ^ {-2\beta} \cdot (2R + 1) ^ {2\beta} \geq s\right\}\right| ds \\
        & \leq \left(\frac{2\beta}{1 - 2\beta}n + 1\right)|\X| + (2R + 1) ^ {2\beta}I_{2\beta}(\X),  \end{aligned}
    \]
    which is indeed finite if  $\beta<1/2$ by \eqref{Ifin} of  Lemma \ref{lem:Minkowski}. Let now $\beta=\alpha-1$ so that $\beta<1/2$ provided $\alpha<3/2$, and then
     $c_1 = \frac{2\beta}{1 - 2\beta}n + 1, c_2 = (2R + 1) ^ {2\beta}$. We can compute the Rényi divergence of $\lambda, \sigma$ to be
    \[
        \begin{aligned}
        D_\alpha(\lambda || \sigma) & = \frac{1}{\alpha - 1}\ln\left(\int_\X \left(\frac{d\lambda}{d\sigma}\right) ^ {\alpha - 1} d\lambda\right) \\
        & \leq \frac{1}{\beta}\ln\left(\int_\X \left(\frac{1}{m\sigma'}\right) ^ {\beta} M ^ \alpha dx\right) \\
        & \leq \frac{1}{\beta} \ln\left(\left(c_1 |\X| + c_2I_{2\beta}(\X)\right)M ^ \alpha m ^ {-\beta}\right) \\
        & = \frac{1}{\beta} \ln\left(c_1|\X| + c_2I_{2\beta}(\X)\right) + \frac{\alpha}{\beta}\ln M - \ln m <\infty.
        \end{aligned}
    \]

    Let now $\mu$ be arbitrary. We will approximate $\mu$ using discrete measures, and we need to show that the transport density will also be approximated. For this, we use the following lemma.
    \begin{lemma}
    \label{lem:convergence}
        Let $\lambda, \mu_k, \mu$ be bounded measures, with $\mu_k$ weakly converging to $\mu$, and let $\sigma_k, \sigma$ be the transport densities associated with $\lambda$ and $\mu_k, \mu$, respectively. Then { $\sigma_k$ weakly converges to $\sigma$}.
    \end{lemma}
    To see how the lemma implies our claim, approximate $\mu$ with discrete measures $\mu_k$, and let $\sigma_k$ be the transport density associated with $\lambda, \mu_k$. Then by the lemma, up to a subsequence, $\sigma_k$ converges weakly to $\sigma$. Since the Rényi divergence is lower semicontinuous in the weak topology \cite[Theorem 19]{Van2014}, we get that
    \[
        D_\alpha(\lambda || \sigma) \leq \liminf_{k \rightarrow \infty} D_\alpha (\lambda || \sigma_k) \leq C_{\X, \alpha, R} + \frac{\alpha}{\alpha - 1} \ln M - \ln m.
    \]
		This concludes the proof when $n \geq 2$. 
		
    Finally, for $n = 1$, let $\mathcal L$ be the Lebesgue measure on $[0, 1]$, and let $\lambda' = \lambda \otimes \mathcal L, \mu' = \mu \otimes \mathcal L$, and $\sigma' = \sigma \otimes \mathcal L$. Then it is not hard to see that $\sigma'$ is the transport density associated with $\lambda', \mu'$, and that $D_\alpha(\lambda' || \sigma') = D_\alpha(\lambda || \sigma)$, which proves Lemma \ref{lem:Renyi} in that case. 
\end{proof}
To finish the proof of Lemma \ref{lem:Renyi}, we need to prove Lemma \ref{lem:convergence}.
\begin{proof}[Proof of Lemma \ref{lem:convergence}]
    {We proceed by contradiction and assume that convergence does not hold. By Prokhorov's theorem, we can extract a subsequence of $\sigma_k$ weakly converging to some measure $\sigma' \neq \sigma$. To simplify notation, we keep denoting this subsequence $\sigma_k$.} Let $\pi_k$ be optimal transport plans from $\lambda$ to $\mu_k$; taking a subsequence again, we can assume $\pi_k$ weakly converges to some transport plan $\pi$ between $\lambda$ and $\mu$, which is optimal, according to \cite[Theorem 1.50]{Santam}. The transport density can be described as a functional acting on a compactly supported continuous function $\varphi$, with value equal to
    \[
        \langle \sigma_k, \varphi\rangle = \int_{\R ^ n \times \R ^ n}\int_0 ^ 1 |x - y| \varphi(ty + (1 - t)x) dt d\pi_k.
    \]
    Note that the integrated function is clearly continuous and bounded, so we can simply compute:
    \[
        \begin{aligned}
        \langle \sigma', \varphi\rangle & = \lim_{k \rightarrow \infty} \int_{\R ^ n \times \R ^ n} \int_0 ^ 1 |x - y| \varphi(ty + (1 - t)x) dt d\pi_k \\
        & = \int_{\R ^ n \times \R ^ n} \int_0 ^ 1 |x - y|\varphi(ty + (1 - t)x) dt d\pi = \langle \sigma, \varphi\rangle,
        \end{aligned}
    \]
    which proves that $\sigma = \sigma'$. This {gives the desired contradiction}.
\end{proof}
\begin{remark}\label{rem:support}
 Let us illustrate the difficulty of considering measures $\lambda$ and $\mu$ such that $\lambda\perp\mu$ but $\spt\lambda\cap\spt \mu\neq \emptyset$. Take  $\lambda=\chi_{(0,1)\times(0,1)}$ and $\mu=\sum_{k\in \N} 2^{-k} \delta_{y_k}$ where $y_k$ are dense in $(0,1)\times (0,1)$. In this case it would be natural to consider in the proof above, $\X= (0,1)\times(0,1)\backslash \spt \mu$ however $\X\backslash \X_r=\X$ for every $r>0$ so that no estimate like \eqref{difvolrectif} can hold.
\end{remark}

Using the upper bound on $D_\alpha(\lambda || \sigma)$, we can transfer  the  bound \eqref{eq:stabsigma} involving integration with respect to  $\sigma$ to a bound involving integration with respect to $\lambda$.
\begin{lemma}
\label{lem1}
    Under the assumptions of Theorem \ref{thm1}, for every $p>3$ there exists a constant $C=C(R,p)>0$ such that for every $f \in L ^ 1(\lambda)$,
\[
    \|f\|_{L ^ 1(\lambda)} \leq C \left( |\X| + I_{2/(p - 1)}(\X)\right) ^ {1 - 1 / p}\frac{M}{m ^ {1 / p}}\|f\|_{L ^ p(\sigma)}.
\]
\end{lemma}
\begin{proof}
    Let  $p>3$ and $p'$ be the Hölder conjugate of $p$, so that $1 / p + 1 / p' = 1$. Notice in particular that $p' < 3 / 2$. By Hölder inequality, we see that
    \[
    \begin{aligned}
    \|f\|_{L ^ 1(\lambda)} & = \left\|f \cdot \frac{d\lambda}{d\sigma}\right\|_{L ^ 1(\sigma)} \leq \|f\|_{L ^ p(\sigma)} \left\|\frac{d\lambda}{d\sigma}\right\|_{L^{p'}(\sigma)}
    \\
    & = \|f\|_{L ^ p(\sigma)} \left(\int_\X \left(\frac{d\lambda}{d\sigma}\right) ^ {p' - 1} d\lambda\right) ^ {1 / p'}
    = \|f\|_{L ^ p(\sigma)} \left(e ^ {(p' - 1)D_{p'}(\lambda || \sigma)}\right) ^ {1 / p'} \\
    & \stackrel{\eqref{eq:Renyi}}{\les_{R,p}}  \left( |\X| + I_{2/(p - 1)}(\X)\right) ^ {1 - 1 / p} \frac{M}{m ^ {1 / p}}\|f\|_{L ^ p(\sigma)}.
    \end{aligned}
    \]
\end{proof}
Theorem \ref{thm1} is now a simple corollary of the lemma.
\begin{proof}[Proof of Theorem \ref{thm1}]
    For $\alpha>3$ we apply Lemma \ref{lem1} to $f = |\nabla \psi - \nabla \phi|\le 2$ and $p=\alpha$ to get that
    \[
        \begin{aligned}
        \|\nabla \psi - \nabla \phi\|_{L ^ 1(\lambda)} & \leq C \|\nabla \psi - \nabla \phi\|_{L ^ p(\sigma)} \leq C \|\nabla \psi - \nabla \phi\|_{L ^ 2(\sigma)} ^ {2 / p} \\
        & \stackrel{\eqref{eq:stabsigma}}{\leq} C \left(\int_{\R^n}(\psi-\phi)d(\mu-\lambda)\right) ^ {1 / p},
        \end{aligned}
    \]
    for a constant $C>0$ depending on $\alpha, R, \X, m$ and $M$.
\end{proof}
\begin{remark}\label{optimality}
 Let $n=1$. Set  $\lambda=\chi_{(0,1)}$, $\mu=\chi_{(1,2)}$ and $\psi(x)=x$ which is an optimal Kantorovich potential for $W_1(\lambda,\mu)$. For $\eps\ll1$ set
 \[
    \phi(x)    =\begin{cases}
                -x & \textrm{if } x\le \eps\\
                x - 2\eps & \textrm{if } x\ge \eps.
            \end{cases}
 \]
We then have
\[
    \|\psi'-\phi'\|_{L^1(\lambda)}=2 \eps \qquad \textrm{and} \qquad  W_1(\lambda,\mu)-\int_{0}^2 \phi d(\mu-\lambda)=2\eps^2-2\int_0^{\eps} xdx=\eps^2.
\]
Therefore \eqref{eq:thm1} cannot hold with $\alpha >1/2$. This example may easily be extended to arbitrary space dimension.
\end{remark}

We now prove that as in \cite{LetMer} for the case $p=1$, stability still holds for some measures which may not be bounded from below or have unbounded support. This comes however with a worse stability exponent. For $\lambda$ absolutely continuous and $R>0$ we introduce the following conditions:
\begin{enumerate}[(H1)]
\item\label{H1} $\X=\spt \lambda\subset B_R$ and $\lambda$ is log-concave;
\item\label{H2} $\lambda$ is $\kappa$-uniformly log-concave;
\item\label{H3}  $\X=\spt \lambda\subset B_R$, $\partial \X$ is $C^2$ and there exist constants $0 < m, M < \infty$, and $\delta > 0$ such that $md(x, \partial \X) ^ \delta \leq \lambda \leq M$ on $\X$.
\end{enumerate}

    Under assumptions (H1) or (H2), we write $\lambda=\exp(-V)$ for some convex function $V$ and assume without loss of generality that the minimum of $V$ is attained at $x=0$. We then fix $r_0>0$ such that $B_{2r_0}\subset  \X$, set $m_0=\inf_{B_{r_0}} \lambda$ and $M=\exp(-V(0))=\max_\X \lambda$.
\begin{theorem}
    \label{thm2}
   Assume that  $\lambda, \mu$ are probability measures with disjoint supports, that there exists $R>0$ such that $\spt \mu \subset B_R$ and that moreover one of the conditions (H1), (H2) or (H3) holds. Then, for every $\alpha>4$ for (H1) and (H2) or $\alpha>3+\delta$ for (H3), there exists a constant $C>0$ with the following property. If $\psi$ is an optimal Kantorovich potential associated with $W_1(\lambda, \mu)$, and $\phi$ is any $1$-Lipschitz function, then
    \begin{equation}\label{stab1unbounded}
        \|\nabla \psi - \nabla \phi\|^\alpha_{L ^ 1(\lambda)} \leq C \int_{\R^n}(\psi-\phi)d(\mu-\lambda).
    \end{equation}
     The constant $C$ depends on $R$, $m_0$, $M$ and $\alpha$ for (H1) and (H2) as well as $\kappa$ for (H2) and on $R$, $\X$, $m$, $M$, $\delta$ and $\alpha$ for (H3).
\end{theorem}
\begin{remark}
 Arguing as in \cite{DelalMeri} it is likely that stability may also be obtained for measures $\mu$ with unbounded support provided we assume finiteness of high enough moments.
\end{remark}

% \begin{theorem}
%     \label{thm3}
%     Suppose $\lambda, \mu$ are probability measures with disjoint supports. Suppose in addition $\mu$ has support contained in $X = B_R(0)$, and $\lambda$ is $\kappa$-uniformly log-concave. Let $\psi$ be the Kantorovich potential associated with $(\lambda, \mu)$, and $\phi$ be any $1$-Lipschitz function. Then
%     \[
%         \|\nabla \psi - \nabla v\|_{L ^ 1(\lambda)} \leq C_{\lambda, \alpha, R} \left(\int_{\R^n}(\psi-\phi)d(\mu-\lambda) \right) ^ {\alpha},
%     \]
%     for every $0 < \alpha < 1 / 4$.
% \end{theorem}
% \begin{theorem}
%     \label{thm4}
%     Suppose $\lambda, \mu$ are probability measures with disjoint supports contained in $X = B_R(0)$, and in addition $\lambda$ is absolutely continuous, supported on a set $\X$ with $C ^ 2$ boundary, and there exist constants $0 < m, M < \infty$, and $\delta > 0$ such that $md(x, \partial \Omega) ^ \delta \leq \lambda \leq M$ on $\Omega$. Let $\psi$ be the Kantorovich potential associated with $(\lambda, \mu)$, and $\phi$ be any $1$-Lipschitz function. Then
%     \[
%         \|\nabla u - \nabla v\|_{L ^ 1(\lambda)} \leq C_{\lambda, \alpha, R} \left(\int_{\R^n}(\psi-\phi)d(\mu-\lambda) \right) ^ {\alpha},
%     \]
%     for every $0 < \alpha < 1 / 4$.
% \end{theorem}
Our method of proving these theorems will involve restricting the source measure $\lambda$ to "nicer" domains, on which Theorem \ref{thm1} holds. We first prove the following lemma.
\begin{lemma}
\label{lem3}
    Let $\lambda, \mu$ be two probability measures, with $\lambda$ absolutely continuous. Let $T$ be an optimal transport map from $\lambda$ to $\mu$, and let $\sigma$ be the transport density associated with $W_1(\lambda, \mu)$. For $A \subseteq \mathbb R ^ n$, let $\lambda_A = \lambda\restr A, \mu_A = T \# \lambda_A$, and let $\sigma_A$ be the transport density associated with $W_1(\lambda_A, \mu_A)$. Then we have $\sigma_A \le \sigma$.
\end{lemma}
\begin{proof}
    Notice that $T$ is also an optimal transport map from $\lambda_A$ to $\mu_A$. With obvious notation (recall \eqref{defmuj}) we have for every non-negative test function $\zeta$, and every $t\in [0,1]$,
    \[
     \int_{\R^n} \zeta d|(j_A)_t|=\int_{\R^n}  \zeta (T_t)  |T-x| d\lambda_A\le \int_{\R^n} \zeta (T_t)  |T-x| d\lambda
     =\int_{\R^n} \zeta d|j_t|.
    \]
    By \eqref{eq:transportdens} this concludes the proof.
\end{proof}
\begin{proof}[Proof of Theorem \ref{thm2}]

    We will prove  results similar to Lemma \ref{lem1}, which will then imply the statements thanks to \eqref{eq:stabsigma}. In this proof we write for compactness
    $a\les_C b$ to indicate an estimate which holds up to a constant which depends on the same parameters as the constant $C$ from \eqref{stab1unbounded} i.e. $R$, $m_0$, $M$ and $\alpha$ for (H1) etc... We start by observing that since
    \[
     \|\nabla \psi - \nabla \phi\|_{L ^ 1(\lambda)}\le 2,
    \]
we may assume without loss of generality that the right-hand side of \eqref{stab1unbounded} is small. By \eqref{eq:stabsigma} this means that we can assume that (here as usual $\sigma$ denotes the transport density associated to $W_1(\lambda,\mu)$)
\begin{equation}\label{hypsmall}
 \|\nabla \psi-\nabla \phi\|_{L^2(\sigma)}^2\ll_C 1.
\end{equation}
{\it Step 1.} We first argue for (H1) and (H2) which we treat simultaneously.

    Let $T$ be an optimal transport map for $W_1(\lambda,\mu)$. For $r\ge r_0, m\le m_0$, let $B_{r, m} = \{x \in B_r : \lambda(x) \geq m$\}. Note that since $\lambda$ is log-concave, $B_{r, m}$ is a convex set. We also let $\lambda_{r, m} = \lambda\restr{B_{r, m}} / \lambda(B_{r, m}), \mu_{r, m} = T \# \lambda_{r, m}$, and finally $\sigma_{r, m}$ be the transport density associated to $W_1(\lambda_{r, m}, \mu_{r, m})$. Notice that since $r\ge r_0$ and $m\le m_0$, we have $B_{r_0}\subset B_{r,m}$ and thus
    \begin{equation}\label{lowerlambda}
     \lambda(B_{r,m})\ge \lambda(B_{r_0})\ge m_0|B_{r_0}|\ges_C 1.
    \end{equation}
    Let $p= \alpha-1>3$. We claim that for every
    $f\in L ^ \infty(\lambda)$, such that $\|f\|_{\infty} \leq 2$ and $\|f\|_{L^p(\sigma)}\ll_C1$, we have
    \begin{equation}\label{claimstabunbounded}
{\|f\|_{L^1(\lambda)}}\les_C \|f\|_{L^p(\sigma)}^{p/(p+1)}.
    \end{equation}
 To this aim, we write
    \be
    \label{thmeq1}
        \|f\|_{L ^ 1(\lambda)} = \|f\|_{L ^ 1(\lambda, B_{r, m})} + \|f\|_{L ^ 1(\lambda, B_r \backslash B_{r, m})} + \|f\|_{L ^ 1(\lambda, \mathbb R ^ d \backslash B_r)}.
    \ee
    We estimate the three terms separately.
      For the first term, we use Lemmas \ref{lem1} and \ref{lem3} as well as \eqref{lowerlambda} to get
    \begin{multline*}
            \|f\|_{L ^ 1(\lambda, B_{r, m})} \leq \|f\|_{L ^ 1(\lambda_{r, m})} \les_C (|B_{r,m}|+ I_{2/(p-1)}(B_{r,m}))^{1-\frac{1}{p}} m ^ {-1 / p} \|f\|_{L ^ p(\sigma_{r, m})} \\
            \les_C (r^n+ I_{2/(p-1)}(B_{r,m}))^{1-\frac{1}{p}} m ^ {-1 / p} \|f\|_{L ^ p(\sigma)}.
    \end{multline*}
     Since $B_{r,m}\subset B_r$ and $B_{r,m}$ is convex, we have by \eqref{MinkowskiConvex} that for $r\ge r_0$,
     \[
    I_{2/(p-1)}(B_{r,m})\les_p r^{n-1}\le  r_0^{-1} r^n\les_C r^n.
     \]
Therefore, we get
    \begin{equation}\label{firsttermthmeq1}
        \|f\|_{L ^ 1(\lambda, B_{r, m})} \les_C m ^ {-1 / p} r ^ n \|f\|_{L ^ p(\sigma)}.
    \end{equation}
For the  second term in \eqref{thmeq1}, we simply observe that since $B_{r,m}\subset B_r$,
    \begin{equation}\label{secondtermthmeq1}
     \|f\|_{L ^ 1(\lambda, B_r \backslash B_{r, m})}\les m r^n.
    \end{equation}
 We finally turn to the last term in \eqref{thmeq1}. If we work under hypothesis (H1) we simply choose $r=R\ge r_0$, for which this term drops. If instead we assume (H2), we observe that since $V$ is minimal at $0$ and uniformly $\kappa-$convex, we have
 \[
  V(x)\ge V(0)+ \frac{\kappa}{2} |x|^2
 \]
 and thus
 \[
  \lambda(x)\le M \exp\lt(-\frac{\kappa}{2} |x|^2\rt).
 \]
Therefore
\begin{equation}\label{thirdtermthmeq1}
 \|f\|_{L ^ 1(\lambda, \mathbb R ^ d \backslash B_r)}\les \lambda(\mathbb R ^ d \backslash B_r)\les_C r^{n-2} \exp\lt(-\frac{\kappa}{2} r^2\rt).
\end{equation}
Combining \eqref{firsttermthmeq1}, \eqref{secondtermthmeq1} and \eqref{thirdtermthmeq1} together with \eqref{thmeq1} we find (recall that under (H1) the third term drops) for $r\ge r_0$ and $m\le m_0$,
    \[
    \|f\|_{L ^ 1(\lambda)} \les_C \left(\frac{r ^ n\|f\|_{L ^ p(\sigma)}}{m ^ {1 / p}} + mr ^ n + r ^ {n - 2} \exp\lt(-\frac{\kappa}{2} r^2\rt)\right).
    \]
    We now conclude the proof of \eqref{claimstabunbounded} under hypothesis (H2) since the argument for (H1) is simpler. We first optimize in $m$ by choosing $m = \|f\|_{L ^ p(\sigma)} ^ {p / (p + 1)}$. Notice that thanks to the hypothesis $\|f\|_{L^p(\sigma)}\ll1$ we have $m\le m_0$. We get
    \[
    \|f\|_{L ^ 1(\lambda)} \les_C \left(r ^ n \|f\|_{L ^ p(\sigma)} ^ {p / (p + 1)} + r ^ {n - 2} \exp\lt(-\frac{\kappa}{2} r^2\rt)\right).
    \]
    We then optimize in $r$ by setting $r = (2\kappa ^ {-1}|\log\|f\|_{L ^ p(\sigma)}|) ^ {1 / 2n}$. Again $r\ge r_0$ thanks to the hypothesis $\|f\|_{L^p(\sigma)}\ll1$. We finally get
    \[
    \|f\|_{L ^ 1(\lambda)} \les_C \|f\|_{L ^ p(\sigma)} ^ {p / (p + 1)} |\log \|f\|_{L ^ p(\sigma)}| ^ {n / 2}.
    \]
   Up to changing slightly the value of $p$, this concludes the proof of \eqref{claimstabunbounded}.\\
   Finally, using \eqref{claimstabunbounded} with  $f = \nabla \psi - \nabla \phi$ which satisfies  $\|f\|_{L^p(\sigma)}\ll1$ by \eqref{hypsmall}, in combination with \eqref{eq:stabsigma} concludes the proof of \eqref{stab1unbounded} under hypothesis (H1) or (H2) {since we get
\begin{align*}
||f||_{L^1(\lambda)} &\les \|f\|_{L ^ p(\sigma)} ^ {p / (p + 1)} \\
&\les ||f||_{L ^ 2(\sigma)} ^ {2 / (p + 1)} \\
&\les \left(\int_{\R^n}{(\psi - \phi) d(\mu -\lambda)}\right)^{1/(p+1)}
\end{align*}
and $p+1 = \alpha$.}

   {\it Step 2.} We now assume (H3).
%     By the coarea formula we can conclude that the integral
%     \[
%         I_\alpha(\X_t) = \int_{\X_t} d(x, \partial \X_t) ^ {-\alpha} dx,
%     \]
%     stays bounded as $t \rightarrow 0$, for a fixed $0 \leq \alpha < 1$. Indeed, if we let $f = d(x, \partial \X)$, then $|\nabla f| = 1$; we denote the levels sets of $f$ as $\X_{(=t)}$. If $t, s$ are small enough, say smaller than $t_0$, then $(\X_t)_{(=s)} = \X_{(=(t + s))}$. The coarea formula gives us
%     \[
%         \begin{aligned}
%         \int_{\X_t \backslash \X_{t_0}} d(x, \partial X_t) ^ {-\alpha} dx & = \int_0 ^ {t_0} \int_{f ^ {-1}(s)} d(x, \partial \X_t) ^ {-\alpha} d\mathcal H_{n - 1}(x) ds \\
%         & = \int_0 ^ {t_0} \int_{(\X_t)_{(=s)}} d(x, \partial \X_t) ^ {-\alpha} d\mathcal H_{n - 1}(x) ds \\
%         & = \int_0 ^ {t_0} \int_{(\X_t)_{(=s)}} s ^ {-\alpha} d\mathcal H_{n - 1}(x) ds = \int_0 ^ {t_0} s ^ {-\alpha} \mathcal H_{n - 1}\left((\X_t)_{(=s)}\right) ds \\
%         & = \int_0 ^ {t_0} s ^ {-\alpha} \mathcal H_{n - 1}\left(\X_{(=(t + s))}\right) ds,
%         \end{aligned}
%     \]
%     and the last integrand stays bounded as $t \rightarrow 0$, since $\alpha < 1$ and $\mathcal H_{n - 1}(\X_{(=t)})$ stays bounded when $t \rightarrow 0$. On $\X_{t_0}$, $d(x, \partial \X_t) ^ {-\alpha}$ is  clearly bounded by $t_0 ^ {-\alpha}$, so we conclude that $I_\alpha(\X_t)$ stays bounded as $t \rightarrow 0$.
%
    We now set $p=\alpha-\delta>3$ and claim that for every $f\in L^\infty(\lambda)$ with $\|f\|_\infty\le 2$ and $\|f\|_{L^p(\sigma)}\ll1$, we have
    \begin{equation}\label{claimH3}
     \|f\|_{L ^ 1(\lambda)}\les_C \|f\|_{L ^ p(\sigma)}^{\frac{p}{p+\delta}}.
    \end{equation}
For this we recall the definition \eqref{Xr} of $\X_r$ and decompose for $r>0$
    \[
      \|f\|_{L ^ 1(\lambda)} = \|f\|_{L ^ 1(\lambda, \X_r)} + \|f\|_{L ^ 1(\lambda, \X \backslash \X_r)}.
    \]
Since $\lambda\ge m r^{-\delta}$ on $\X_r$, if $r$ is small enough so that $\lambda(\X_r)\ges 1$, we have arguing as in Step 1 by combining Lemma \ref{lem3} and Lemma \ref{lem1},
\[
 \|f\|_{L ^ 1(\lambda, \X_r)}\les_C (|\X_r|+ I_{2/(p-1)}(\X_r))^{1-1/p} r^{-\delta/p}\|f\|_{L ^ p(\sigma)}\stackrel{\eqref{limsupI}}{\les_C}r^{-\delta/p}\|f\|_{L ^ p(\sigma)}.
\]
Since
\[
 \|f\|_{L ^ 1(\lambda, \X \backslash \X_r)}\les_C |\X \backslash \X_r|\stackrel{\eqref{difvolrectif}}{\les_C} r,
\]
we find for $r$ small enough,
\[
 \|f\|_{L ^ 1(\lambda)}\les_C r^{-\delta/p}\|f\|_{L ^ p(\sigma)}+ r.
\]
Optimizing in $r$ by choosing $r=\|f\|_{L ^ p(\sigma)}^{\frac{p}{p+\delta}}\ll1$ by the hypothesis $\|f\|_{L^p(\sigma)}\ll1$, we get \eqref{claimH3}.
The proof of \eqref{stab1unbounded} is then concluded as above by combining this with \eqref{eq:stabsigma}.\\
\end{proof}
\section{Convolution inequality}\label{sec:convol}
In this section we turn to the study of \eqref{eq:contraction}. We recall the definition \eqref{def:deltaepsmunu} of $\delta_\eps(\lambda,\mu)$. Also recall that for a measure $\mu$ and $z\in \R^n$, we define the translated measure $\mu^z$ by $\mu^z(A)=\mu(A-z)$ for $A\subset \R^n$. In particular, we have
\[
 \mu_\eps=\int_{\R^n} \mu^z d\rho_\eps(z).
\]
\begin{proof}[Proof of Theorem \ref{theo:equal:intro}]
As a warm-up we first prove \eqref{eq:contraction} in full generality, i.e. without any assumptions on $\lambda$ (in particular it does not need to be absolutely continuous) and $\mu$ or $\rho$. Although it is well-known, we provide a proof based on the dual problem since it will be the starting point of both the rigidity statement from Theorem \ref{theo:equal:intro} as well as  its quantitative counterparts.\\ 
Since $W^p_p(\lambda^z,\mu^z)=W^p_p(\lambda,\mu)$ for every $z\in \R^n$, we have for every $\phi$
\begin{multline*}
    W^p_p(\lambda,\mu)=\int_{\R^n} W^p(\lambda^z,\mu^z) d\rho_\eps(z)\ge \int_{\R^n} \lt[\int_{\R^n} \phi^c d\lambda^z + \int_{\R^n} \phi d\mu^z\rt] d\rho_\eps(z)\\
    =\int_{\R^n} \phi^c d\lambda_\eps+ \int_{\R^n} \phi d\mu_\eps.
\end{multline*}
Taking the supremum over $\phi$ concludes the proof of \eqref{eq:contraction}.\\
Assume now that $\lambda$ is absolutely continuous, and that either $\spt \rho=\R^n$, or $\X$ is connected and moreover that there exist optimal Kantorovich potentials  $\psi$ for $W_p(\lambda,\mu)$ and $\psi_\eps$ for $W_p(\lambda_\eps,\mu_\eps)$.
In particular  we have
\[
 \int_{\R^n} \psi_\eps^c d\lambda_\eps+ \int_{\R^n} \psi_\eps d\mu_\eps=W_p^p(\lambda_\eps,\mu_\eps).
\]
{To prove Theorem \ref{theo:equal:intro}, we focus on the only if part, since the if part is immediate by a direct coupling argument.}
If $\delta_\eps(\lambda,\mu)=0$, then we have equality in the above inequality for the choice $\phi=\psi_\eps$, and thus $\psi_\eps$ is an optimal Kantorovich potential for $W_p(\lambda^z,\mu^z)$ for every $z\in \spt \rho_\eps$. If $\psi$ is an optimal Kantorovich potential for $W_p(\lambda,\mu)$, we claim that $\nabla \psi^c$ is constant in $\X$. To prove this claim we notice that for every $z\in \R^n$, $\psi^z$ is also optimal for $W^p_p(\lambda^z,\mu^z)$. We thus have for a.e. $x\in \spt \lambda^z$ (for $p=1$ the equality holds for $x\in \spt \sigma^z$ but $\spt\lambda^z\subset \spt \sigma^z$ by Lemma \ref{lem:qualitabscont}),
\[
  \nabla \psi_\eps^c(x)= \nabla(\psi^z)^c(x) =\nabla \psi^c(x-z).
\]
Making the change of variable $y=x-z$, this is equivalent to
\[
 \nabla \psi^c(y)=\nabla \psi_\eps^c(y+z),
\]
for $y\in \spt \lambda$ and $z\in \spt \rho_\eps$. If $\spt \rho_\eps=\R^n$ this concludes the proof of the claim. Assume now that $\X$ is connected.  Letting $f(x)=\nabla \psi_\eps^c(x)$ and $\xi(y)=\nabla \psi^c(y)$, since $\inf_{B_1} \rho>0$, we have
 for a.e. $x\in \R^n$,
\[
 \xi =f(x) \qquad \textrm{a.e. in } B_\eps(x)\cap \X.
\]
If $x,x'$ are Lebesgue points of $f$ such that $|B_\eps(x)\cap B_\eps(x')\cap \X|>0$ we thus have
\[
 f(x')=\xi =f(x) \qquad \textrm{a.e. in } B_\eps(x)\cap B_\eps(x')\cap \X,
\]
so that $f(x)=f(x')$. We conclude that $\xi$ is constant in $(B_\eps(x)\cup B_\eps(x'))\cap \X$. Taking for $y,y'\in \X$ a chain of consecutively overlapping balls of radius $\eps$ joining $y$ to $y'$ we get that $\xi$ is constant on $\X$ as claimed.\\
We now argue separately in the cases $p>1$ and $p=1$. In the first case, since $\nabla \psi^c$ is constant in $\X$ also (recall the definition \eqref{def:Phip} of $\Phi_p$)
\[
 y-T(y)=\Phi_p(\nabla \psi^c(y)),
\]
is constant. Letting $-z$ be this constant we have  $T\# \lambda=\lambda^{z}$ which concludes the proof in this case.\\
If instead $p=1$, from the fact that $\nabla \psi^c=-\nabla \psi$ is constant we deduce that up to an additive constant, $\psi(x)=\ang{x, e}$ for some $e\in \bS^{n-1}$. Let $\vhi$ be monotone in the direction $e$. We then have
\[
 \int_{\R^n} \vhi d\mu -\int_{\R^n} \vhi d\lambda=\int_{\R^n} \vhi(T(x))-\vhi(x) d\lambda\ge 0,
\]
where we used that $T(x)-x= |T(x)-x| e$ and thus $\vhi(T(x))-\vhi(x)\ge 0 $ by monotonicity of $\vhi$.
\end{proof}
We now turn to the quantitative version of \eqref{eq:contraction}. Let us first connect it with the quantity $\Lambda_\eps(\xi,f)$ (recall \eqref{defLambda}) from Section \ref{sec:2point}. In order to have a unified presentation for $p>1$ and $p=1$ let us set
\begin{equation*}\label{hatphip}
 \hat{\Phi}_p(z)=\begin{cases}
                  \Phi_p(z) & \textrm{if } p>1\\
                  z & \textrm{if } p=1.
                 \end{cases}
\end{equation*}
Recall the definition \eqref{def:Kantorovichfunct} of $F_{\lambda,\mu}$. Letting  $\psi$, respectively $\psi_\eps$, be a maximizer for $F_{\lambda,\mu}$, respectively $F_{\lambda_\eps,\mu_\eps}$; setting $\phi= \psi^c$ and $\phi_\eps=\psi_\eps^c$, we make the following assumption:
\begin{hypothesis}\label{hypstrongconvex}
There exists $ C>0$ and $\alpha\ge 1$ such that for every $z\in \spt \rho_\eps$,
\begin{equation}\label{strongconvex}
 F_{\lambda^z,\mu^z}(\psi^z)-F_{\lambda^z,\mu^z}(\psi_\eps)\ge  C \lt(\int_{\R^n}|\hat{\Phi}_p(\nabla \phi^z)-\hat{\Phi}_p( \nabla \phi_\eps)|^p d\lambda^z\rt)^\alpha.
\end{equation}
\end{hypothesis}
\begin{proposition}\label{connectLambda}
 If Hypothesis \ref{hypstrongconvex} holds then
 \begin{equation}\label{Lambdadelta}
  \delta_\eps^{\frac{1}{\alpha}}(\lambda,\mu)\ge C^{\frac{1}{\alpha}} \Lambda_\eps(\hat{\Phi}_p(\nabla \psi^c),\hat{\Phi}_p(\nabla \psi^c_\eps)).
 \end{equation}
\end{proposition}
\begin{proof}
 We start by observing that if $\psi$ is optimal for $F_{\lambda,\mu}$ then for $z\in \spt \rho_\eps$ and $x\in \spt \lambda^z$ we have
 \begin{equation}\label{equalpsizc}
  (\psi^c)^z(x)=(\psi^z)^c(x).
 \end{equation}
Indeed, writing that $x=z+ \hat{x}$ for some $\hat{x}\in \spt \lambda$, we have
\begin{multline*}
 (\psi^c)^z(x)=\psi^c(x-z)=\psi^c(\hat{x})=\inf_{y} \frac{1}{p}|\hat{x}-y|^p -\psi(y)\\
 =\inf_{y} \frac{1}{p}|\hat{x}+z-y|^p -\psi(y-z)=(\psi^z)^c(x).
\end{multline*}
We write

\[ \delta_\eps(\lambda,\mu)=W^p_p(\lambda,\mu)- W^p_p(\lambda_\eps,\mu_\eps)=\lt[\int_{\R^n} \psi d\mu + \int_{\R^n} \psi^c d\lambda\rt] -\lt[\int_{\R^n} \psi_\eps d\mu_\eps+ \int_{\R^n}\psi_\eps^c d\lambda_\eps\rt].\]
Since for every $z\in \spt \rho_\eps$,
\[
 \int_{\R^n} \psi d\mu + \int_{\R^n} \psi^c d\lambda=\int_{\R^n} \psi^z d\mu^z + \int_{\R^n} (\psi^c)^z d\lambda^z,
\]
using \eqref{equalpsizc} and $\int_{\R^n} d\rho_\eps(z)=1$ we get
\begin{multline*}
 \int_{\R^n} \psi d\mu + \int_{\R^n} \psi^c d\lambda=\int_{\R^n} \lt[\int_{\R^n} \psi^z d\mu^z + \int_{\R^n} (\psi^z)^c d\lambda^z\rt] d\rho_\eps(z)
 =\int_{\R^n} F_{\lambda^z,\mu^z}(\psi^z) d\rho_\eps(z).
\end{multline*}
Combining this with
\[
 \int_{\R^n} \psi_\eps d\mu_\eps+ \int_{\R^n}\psi_\eps^c d\lambda_\eps=\int_{\R^n}\lt[\int_{\R^n} \psi_\eps d\mu^z+ \int_{\R^n} \psi_\eps^c d\lambda^z\rt] d\rho_\eps(z)=\int_{\R^n}F_{\lambda^z,\mu^z}(\psi_\eps) d\rho_\eps(z),
\]
we find
\[
  \delta_\eps(\lambda,\mu)= \int_{\R^n} \lt[F_{\lambda^z,\mu^z}(\psi^z)-F_{\lambda^z,\mu^z}(\psi_\eps)\rt] d\rho_\eps(z).
\]
Using \eqref{strongconvex} and Jensen's inequality we obtain
\begin{align*}
  \frac{\delta_\eps(\lambda,\mu)}{  C}&\ge \int_{\R^n} \lt(\int_{\R^n}|\hat{\Phi}_p(\nabla \phi^z)-\hat{\Phi}_p( \nabla \phi_\eps)|^p d\lambda^z\rt)^\alpha d\rho_\eps(z)\\
 &\ge\lt(\int_{\R^n} \int_{\R^n}|\hat{\Phi}_p(\nabla \phi^z)-\hat{\Phi}_p( \nabla \phi_\eps)|^p d\lambda^z d\rho_\eps(z)\rt)^\alpha\\
 &\stackrel{x-z=y}{=}\lt(\int_{\R^n} \int_{\R^n}|\hat{\Phi}_p(\nabla \phi(y))-\hat{\Phi}_p( \nabla \phi_\eps(x))|^p d\lambda(y) d\rho_\eps(x-z)\rt)^\alpha\\
 &=\Lambda_\eps(\hat{\Phi}_p(\nabla \phi),\hat{\Phi}_p( \nabla \phi_\eps))^\alpha.
\end{align*}
This concludes the proof of \eqref{Lambdadelta}.
\end{proof}
From \eqref{Lambdadelta} we see that we can obtain quantitative versions of Theorem \ref{theo:equal:intro}, provided we can check Hypothesis \ref{hypstrongconvex} and we have a lower bound on $\Lambda_\eps(\hat{\Phi}_p(\nabla \psi^c),\hat{\Phi}_p( \nabla \psi^c_\eps))$.
Combining Theorem \ref{twopointcompact} and Proposition \ref{connectLambda} we get the following conditional result.
\begin{theorem}\label{theo:conditional}
 Assume that Hypothesis \ref{hyp:lambda} and Hypothesis \ref{hypstrongconvex} hold {(with $r=\eps$ for Hypothesis \ref{hyp:lambda})}. Setting
 \[
  z=\int_{\R^n} \hat{\Phi}_p(\nabla \psi^c) d\lambda,
 \]
we have
 \begin{equation}\label{quantitativeconditional}
  \int_{\R^n} |\hat{\Phi}_p(\nabla \psi^c)-z|^p d\lambda \les M_0(\eps)\tau_{p,\lambda}(\eps)\delta_\eps^{\frac{1}{\alpha}}(\lambda,\mu).
 \end{equation}
Here the implicit constant depends on $ C$ and $\alpha$ from Hypothesis \ref{hypstrongconvex}.
 \end{theorem}
\begin{remark}\label{translatcondi}
For $p>1$,  since $x- \hat{\Phi}_p(\nabla \psi^c)+z=T+z$ is a transport map between $\lambda$ and $\mu^z$,  \eqref{quantitativeconditional}  implies
 \[
  W_p^p(\lambda,\mu^z)\les M_0(\eps)\tau_{p,\lambda}(\eps)\delta_\eps^{\frac{1}{\alpha}}(\lambda,\mu).
 \]
 If instead $p=1$, \eqref{quantitativeconditional} reads
 \[
  \int_{\R^n} |\nabla \psi+z|d\lambda\les M_0(\eps)\tau_{1,\lambda}(\eps)\delta_\eps^{\frac{1}{\alpha}}(\lambda,\mu).
 \]
Notice that while a priori $|z|\neq 1$, since $|\nabla \psi|=1$ a.e. on $\X=\spt \lambda$,
\[
 |\nabla \psi+z|\ge \min_{e\in \bS^{n-1}}|z-e|
\]
which after integration yields
\[
 \min_{e\in \bS^{n-1}}|e+z|\le \int_{\R^n} |\nabla \psi+z|d\lambda\les M_0(\eps) \tau_{1,\lambda}(\eps)\delta_\eps^{\frac{1}{\alpha}}(\lambda,\mu).
\]
Thus if $e$ reaches the minimum on the right-hand side, by triangle inequality we have for some $e\in \bS^{n-1}$,
\begin{equation}\label{stab1e}
 \int_{\R^n} |\nabla \psi-e|d\lambda\les M_0(\eps)\tau_{1,\lambda}(\eps)\delta_\eps^{\frac{1}{\alpha}}(\lambda,\mu).
\end{equation}

\end{remark}

We now provide non-conditional results. At the moment we can obtain them only in the case $p=2$ or $p=1$, for which $\hat{\Phi}_p(z)=z$.
%
% \begin{proposition}
%     $(\lambda, \mu)$ achieve equality if and only if for some $e \in \mathbb R ^ n$, we have $
%     \eta = p_e \# \lambda = p_e \# \mu$, and denoting $\lambda = \lambda_y \otimes \eta, \mu = \mu_y \otimes \eta$, $\lambda_y$ is stochastically dominated by $\mu_y$ on $p_e ^ {-1}(y)$ for almost all $y$.
% \end{proposition}

\subsection{Quantitative convolution inequality for $p=2$}\label{subsec:stabp}
We start with the quadratic case $p=2$. In this case we report the following result from  \cite[Corollary 1.31]{delalande2022quantitative} (under the hypothesis that $\X$ is convex) and its extension to John domains in \cite[Theorem 2.8]{FordThesis}, see also the proof of \cite[Corollary 5.6]{delalande2022quantitative}. The proof is very similar to  the proof in \cite{DelalMeri} with a slightly different presentation. To obtain the conclusion in this form we used that both $\psi$ and $ \psi_\eps$ satisfy \eqref{eq:Lippsic}.
\begin{proposition}\label{prop:compactconvexcorrect}
Let $p=2$, and assume that $\X=\spt \lambda$ is  bounded Lipschitz and connected. Assume in addition that there exist constants $0 < m\le M$ such that $m\le \lambda\le M$ on $\X$, and there exists $R>0$ such that $\spt\lambda\cup\spt \mu \subset B_R$. Then there exists $C=C(R, m,M,\X)$ such that for  every maximizer $\psi$ of $F_{\lambda,\mu}$ and every $c-$concave function $\phi\in C(B_R)$,
\[
 F_{\lambda,\mu}(\psi)-F_{\lambda,\mu}(\phi)\ge  C \lt(\int_{\R^n}|\nabla  \psi^c- \nabla \phi^c|^2 d\lambda\rt)^3.
\]
\end{proposition}

\begin{proof}[Proof of Theorem \ref{theo:introsimplequadratic}]
  Let $R>0$ be such that $\spt \lambda_\eps\cup \spt \mu_\eps\subset B_R$. For every $z\in \spt \rho_\eps$ we may apply Proposition \ref{prop:compactconvexcorrect} to $\phi=\psi_\eps^{-z}$ to obtain that Hypothesis \ref{hypstrongconvex} holds with $\alpha=3$. Now since $\X$ is   Lipschitz,  we may apply Lemma \ref{lem:tauplambda} to get $\tau_{2,\lambda}\les_\X \eps^{-(n+1)}$ and $M_0(\eps)\les_\X M$. Combining this with \eqref{quantitativeconditional} concludes the proof.
\end{proof}
We now investigate the case when $\lambda$ is a Gaussian. To guarantee that Hypothesis \ref{hypstrongconvex} holds we will rely on the following result from  \cite[Proposition 1.20]{delalande2022quantitative}, see also \cite{gigli2011holder,hutter2021minimax,merigot2020quantitative,BaMan}.
\begin{proposition}\label{prop:convLipcorrect}
 Assume that the optimal Kantorovich potential $ \psi$ for $W_2(\lambda,\mu)$   is in $C^1(\R^n)$. Then for every $c-$concave function $\phi\in C^{1}(\R^n)$ such that $x-\nabla {\phi}^c$  is $K-$Lipschitz, we have
 \[F_{\lambda,\mu}(\psi)-F_{\lambda,\mu}(\phi)\ge 2 K \int_{\R^n}|\nabla  \psi^c- \nabla \phi^c|^2 d\lambda.
 \]
\end{proposition}
\begin{remark}
Notice that a limiting aspect of this statement is that here  we need to require not that $\nabla \psi^c$ is Lipschitz (which is a statement about the optimal transport  map from $\lambda$ to $\mu$) but that the competitor  $\nabla {\phi}^c$ is Lipschitz. This reduces a lot the applicability of this result. This restriction of the method of proof has been already pointed out, see e.g. \cite[p.8 \& p. 11]{BaMan}. Notice that in the framework of \cite{BaMan}, the problem is somewhat symmetric in the roles of $\psi$ and $\phi$, so they could go around this issue and have a hypothesis only on smoothness of $\psi$, see in particular \cite[Theorem 3]{BaMan} and \cite[p.16]{BaMan}.
\end{remark}
As in the proof of Theorem \ref{theo:introsimplequadratic}, our aim is to apply Proposition \ref{prop:convLipcorrect} to $ \phi=\psi_\eps$ where $\psi_\eps$ is an optimal Kantorovich potential for $W_2(\lambda_\eps,\mu_\eps)$. Notice that if $T_\eps$ is the optimal transport map from $\lambda_\eps$ to $\mu_\eps$ we have $T_\eps=x-\nabla \psi_\eps^c$, recall \eqref{PhipT}. In order to find a framework where we can guarantee that $T_\eps$ is $K-$Lipschitz, we will assume that $\mu$ is $\kappa^{-1}-$log concave and that $\rho_\eps=p_{\sqrt{\eps}}$ is the heat kernel. This allows us to rely on Caffarelli's contraction principle \cite{caffarelli2000monotonicity}, see also \cite{gozlan2025global} for some of the most recent progress on this topic.
\begin{theorem}\label{theo:Cafffarelli}
 Let $\lambda=\exp(-V)$ and $\mu=\exp(-W)$ be two probability densities on $\R^n$. Assume that $\spt \lambda=\R^n$, $\spt \mu$ is convex with non empty interior, that  both $V$ and $W$ are $C^2$ in their domains and
 \begin{equation*}\label{hyp:VW}
  \nabla^2V\le \Sigma^{-1} Id, \qquad \nabla^2W \ge \kappa^{-1} Id,
 \end{equation*}
with $\Sigma,\kappa>0$. Then the optimal transport map $T$ for $W_2(\lambda,\mu)$ is  $\sqrt{ \kappa/\Sigma}-$Lipschitz.
\end{theorem}
We report the following seemingly folklore consequence of the Pr\'ekopa-Leindler inequality.
\begin{lemma}\label{Prek}
 Let $\mu$ be $\kappa^{-1}-$log concave and $p_t$ be the heat kernel then for every $t>0$, $p_t\ast \mu$ is $(\kappa+2t)^{-1}-$log concave.
\end{lemma}
\begin{proof}
 The claim is a direct application of  \cite[Theorem 4.3]{brascamp1976extensions}. We include the computations for the reader's convenience. We write
 \[
  p_t(x)= \frac{1}{Z_t} \exp\lt(-\frac{|x|^2}{4t}\rt) \qquad \textrm{and } \qquad \mu(x)=\exp\lt(-\frac{|x|^2}{2\kappa}- W(x)\rt)
 \]
where $W$ is convex. We then have
\[
 p_t\ast \mu(x)= \frac{1}{Z_t} \int_{\R^n} \exp\lt(-\frac{|x-y|^2}{4t}-\frac{|y|^2}{2\kappa}- W(y)\rt) dy.
\]
Writing that
\[
 \frac{|x-y|^2}{4t}+\frac{|y|^2}{2\kappa}=\frac{1}{2} \frac{|x|^2}{\kappa+2t} + \frac{\kappa+2t}{4tk}\lt|y-\frac{\kappa}{\kappa+2t} x\rt|^2
\]
we obtain
\[
 p_t\ast \mu(x)= \frac{1}{Z_t} \exp\lt(-\frac{1}{2} \frac{|x|^2}{\kappa+2t} \rt) \int_{\R^n} \exp\lt(-\frac{\kappa+2t}{4tk}\lt|y-\frac{\kappa}{\kappa+2t} x\rt|^2- W(y)\rt) dy.
\]
Applying now  \cite[Theorem 4.3]{brascamp1976extensions} with
\[
 \Phi(x,y)=\exp\lt(-\frac{\kappa+2t}{4tk}\lt|y-\frac{\kappa}{\kappa+2t} x\rt|^2\rt) \qquad \textrm{and } \qquad F(x,y)= \exp\lt(-W(y)\rt)
\]
we obtain that (notice that with the notation of \cite{brascamp1976extensions} we have $D=0$)
\[
  \int_{\R^n} \exp\lt(-\frac{\kappa+2t}{4tk}\lt|y-\frac{\kappa}{\kappa+2t} x\rt|^2- W(y)\rt) dy
\]
is log concave and thus $p_t\ast\mu$ is $(\kappa+2t)^{-1}-$log concave.
\end{proof}

We then have
\begin{theorem}\label{theo:stabgauss}
 Let $\lambda$ be a standard Gaussian and $\mu$ be $\kappa^{-1}-$log concave. Assume that $\rho_\eps=p_{\sqrt{\eps}}$ is the heat kernel, then for every $\eps^2\ll \kappa$
 and every $\beta>0$,
\begin{equation}\label{mainthmgauss}
   \min_{z\in \R^n} W^2_2(\lambda,\mu^z)\les_\beta \eps^{-2} \delta_\eps(\lambda,\mu)^{1-\beta}.
  \end{equation}
 \end{theorem}
\begin{proof}
To ease notation, let $t_\eps=\sqrt{\eps}$ and $\Sigma_\eps^2=1+2t_\eps$.  Notice first  that since the heat kernel is a semi-group and since $\lambda=p_{1/2}$ we have
\[\lambda_\eps=p_{t_\eps}\ast p_{1/2}=p_{t_\eps+ 1/2}=\gamma_{\Sigma_\eps}=\exp(-V_\eps),\]
with
\[
 \nabla^2V_\eps=\Sigma_\eps^{-1}.
\]
Moreover, by Lemma \ref{Prek} we get that $\mu_\eps$ is $(\kappa+2t_\eps)^{-1}-$log concave. We now apply a first time Theorem \ref{theo:Cafffarelli} and obtain that the optimal transport map $T=\nabla \psi^c$ for $W_2(\lambda,\mu)$ is $\sqrt{\kappa}-$Lipschitz. Applying a second time Theorem \ref{theo:Cafffarelli} we also get that
 the optimal transport map $T_\eps=\nabla  \psi_\eps^c$ for $W_2(\lambda_\eps,\mu_\eps)$ is $\sqrt{(\kappa+2t_\eps)/\Sigma_\eps}-$Lipschitz. By Proposition \ref{prop:convLipcorrect}, we see that Hypothesis \ref{hypstrongconvex} holds with $\alpha=1$ and
 \[
  C=2K =2\sqrt{(\kappa+2t_\eps)/\Sigma_\eps}\approx \sqrt{\kappa}.
 \]

We may thus apply \eqref{Lambdadelta} and get
\begin{equation}\label{stabfirststepgauss}
 \delta_\eps(\lambda,\mu)\ges \sqrt{\kappa} \Lambda_\eps(\nabla \psi^c, \nabla \psi_\eps^c)=\sqrt{\kappa} \Lambda_\eps(T-x,T_\eps-x).
\end{equation}
By invariance under translation of quadratic optimal transport, we may assume without loss of generality that $T(0)=0$. Since $T$ is $\sqrt{\kappa}-$Lipschitz, we get that $\xi=T-x$ satisfies
\[
 |\xi(x)|\le |x|+|T(x)|\les (1+\sqrt{\kappa})|x|.
\]
% \textcolor{blue}{Mi: I would like now to obtain a bound on $T(0)$ which depends only on $\kappa/\sigma$. On way to obtain a bound is to argue as in e.g. \cite{carlier2024optimal}. We set for $e\in \bS^{n-1}$
% \[
%  K_e=\{x\in \R^n \ : \ang{x,e}\ge \frac{1}{2}|x|\}.
% \]
% Letting $e=T(0)/|T(0)|$ we have for $x\in K_e$ by monotonicity of $T$,
% \[
%  \frac{1}{2}|T(0)||x|\le \ang{T(0),x}\le \ang{T(x),x}\le |T(x)| |x|
% \]
% so that $|T(x)|\ge |T(0)|/2$. We conclude that $T(K_e)\subset B_{|T(0)|/2}^c$ so that
% \begin{equation}\label{estimmonotone}
%  \mu(B_{|T(0)|/2}^c)\ge \lambda(K_e)\ges 1
% \end{equation}
% which gives an implicit bound on $|T(0)|$ (as $\mu(B_r^c)\to0$ as $r\to \infty$). Another line of reasoning takes inspiration from \cite{colombo2021bounds}. Since $p=2$, up to translation we might assume that
% \[
%  z=\int_{\R^n} T-x d\lambda=\int_{\R^n}  yd\mu=0.
% \]
% If $\mu$ is furhter assumed to be isotropic we have
% \[
%  \int_{\R^n} |z|d\mu\le \lt(\int_{\R^n} |z|^2 d\mu\rt)^{\frac{1}{2}}=n^{1/2}.
% \]
% By \cite[Proposition 3.2]{Fathi},
% \[\mu(B_{r+n^{1/2}}^c)\le \exp(- \frac{\kappa^2 r^2}{2}).\]
% Plugging this back into \eqref{estimmonotone} and assuming that $|T(0)|\gg 1$ (otherwise there is nothing to prove) we get
% \[
%  1\les \exp(-\frac{\kappa^2 |T(0)|^2}{8})
% \]
% which finally leads to
% \[
%  |T(0)|\les \kappa.
% \]
% The question is if it is possible to get rid of the isotropy? Does $\kappa^{-1}-$log concavity implies some bound on the co-variance?}
Combining Proposition \ref{prop:twopointgaussian} and \eqref{stabfirststepgauss} we find that there exists $z\in \R^n$ such that for any $\beta\in(0,1)$,
\[
 \int_{\R^n} |T-x-z|^2 d\lambda \les_{\beta,\kappa} \eps^{-2} \delta_\eps^{1-\beta}(\lambda,\mu).
\]
This concludes the proof of \eqref{mainthmgauss}.
\end{proof}

\begin{remark}\label{rem:gaussian case}
 Assume that $\mu=\gamma_\kappa$ is an isotropic and centered Gaussian density. Then  see e.g. \cite[Remark 2.31]{peyre2019computational}
 \[
  W_2^2(\lambda,\mu)= n|1-\kappa|^2.
 \]
In particular if $\rho_\eps= p_{t_\eps}$ is the convolution with the heat kernel at time $t_\eps=\sqrt{\eps}$ we have $\lambda_\eps=\gamma_{\Sigma_\eps}$ and $\mu_\eps=\gamma_{\kappa_\eps}$ with
\[
 \Sigma_\eps^2=1+2t_\eps, \qquad \textrm{and} \qquad \kappa_\eps^2=\kappa^2+ 2t_\eps.
\]
Therefore
\[
 \delta_\eps(\lambda,\mu)=n(|1-\kappa|^2-|\Sigma_\eps-\kappa_\eps|^2)=\lt(1- \frac{|\Sigma_\eps-\kappa_\eps|^2}{|1-\kappa|^2}\rt)W_2^2(\lambda,\mu).
\]
Let
\[
 f(\eps,\kappa)=1- \frac{|\Sigma_\eps-\kappa_\eps|^2}{|1-\kappa|^2}=1- \frac{(1+\kappa)^2}{(\Sigma_\eps +\kappa_\eps)^2}.
\]
We notice first that if $\eps$ is fixed and $\kappa$ tends to $\infty$, then $f\to 0$, so that we cannot hope for a bound of the form \eqref{mainthmgauss} without dependence on $\kappa$. On the other hand, if we assume that $\kappa\le R$ for some $R>0$, we have
\[
 f(\eps,\kappa)\ge \frac{2t_\eps}{R^2+2t_\eps}.
\]
We can thus guess that the optimal prefactor in \eqref{mainthmgauss} should be $ t_\eps^{-1}=\eps^{-\frac{1}{2}}$.

\end{remark}

\subsection{Quantitative convolution inequality for $p=1$}\label{subsec:stapone}
We recall that for $e\in \bS^{n-1}$, we set $p_e(x)=x- \ang{x, e}e$ to be the orthogonal projection on $e^\perp$.
We order each set $p_e ^ {-1}(y)$ by saying that $x_1 \leq x_2$ if and only if $\langle x_1 - x_2, e\rangle \leq 0$.
\begin{proof}[Proof of Theorem \ref{theo:introsimplelinear}]
We write here $a\les_C b$ to indicate an estimate which holds up to constant depending on $\alpha, \X, m, M$ or $R$. By Theorem \ref{thm1}, Hypothesis \ref{hypstrongconvex} holds for any $\alpha>3$. By Lemma \ref{lem:tauplambda} we have that also Hypothesis \ref{hyp:lambda} holds with $\tau_{1,\lambda}(\eps)\les_C \eps^{-n}$ and $M_0(\eps)\les_C 1$ so that by Theorem \ref{theo:conditional}, see also \eqref{stab1e}, there exists $e\in \bS^{n-1}$ such that for $\alpha>3$
\begin{equation}\label{firststep_stab_p1}
 \int_{B_R} |\nabla \psi-e|d\lambda\les_C \eps^{-n} \delta_\eps^{\frac{1}{\alpha}}(\lambda,\mu).
\end{equation}
    We first prove \eqref{quantmonot}. Let $\varphi$ be $1$-Lipschitz and monotone in the direction $e$. Let $T$ be a transport map from $\lambda$ to $\mu$, and let $\tilde T(x) =x + \ang{T(x)-x,e}e$. If $|\nabla \psi(x) - e| \leq \sqrt{2}$, which is equivalent to $\ang{\nabla \psi(x),e}\ge 0$ then  since
    \[T(x) - x= \eta(x)\nabla \psi(x),\]
    for  $\eta(x)=|T(x)-x|>0$, we have that $\tilde T(x) - x$ is a positive multiple of $e$, so $\varphi(\tilde T(x)) \geq \varphi(x)$.
    We may now estimate using $|\nabla \psi(x)|=|e|=1$,
    \begin{align*}
     |T(x)-\tilde T(x)|&=| (T(x)-x) -\ang{T(x)-x,e}e|=\eta(x) |\nabla \psi(x)-\ang{\nabla \psi(x),e}e|\\
     &=|T(x)-x|\sqrt{1-\ang{\nabla \psi(x),e}^2}
     \le |T(x)-x|\sqrt{2(1-\ang{\nabla \psi(x),e}}\\&=|T(x)-x||\nabla \psi(x)-e|,
    \end{align*}
where  we used that $1-\ang{\nabla \psi(x),e}^2=(1-\ang{\nabla \psi(x),e})(1+\ang{\nabla \psi(x),e})\le 2(1-\ang{\nabla \psi(x),e})$.
Since $x,T(x)\in B_R$, this leads to the estimate
\[
 |T(x) - \tilde T(x)|\le 2R|\nabla \psi(x)-e| \qquad \quad \textrm{if }  |\nabla \psi(x)-e|\le \sqrt{2}.
\]
Since $\varphi$ is $1$-Lipschitz, we get when $|\nabla \psi(x) - e| \leq \sqrt{2}$,
    \[
        \varphi(T(x)) \geq \varphi(\tilde T(x)) - 2R|\nabla \psi(x) - e| \geq \varphi(x) - 2R|\nabla \psi(x) - e|.
    \]
     Let $A$ be the set of $x$ for which $|\nabla \psi(x) - e| \ge \sqrt{2}$. By Markov's inequality, we can bound the size of $A$ as
    \[
        \lambda(A) \leq (1/\sqrt{2}) \int_{\R^n}|\nabla \psi(x) - e|d\lambda.
    \]
    On $A$ since $\vhi$ is $1-$Lipschitz, we have  $\varphi(x) - \varphi(T(x))\le 2R$. Therefore, using that $T\#\lambda=\mu$,
    \[
        \begin{aligned}
        \int_{\R^n} \vhi d\lambda-\int_{\R^n}\varphi d\mu & = \int_{\R^n}(\vhi(x)-\varphi(T(x))) d\lambda\\
        &\leq 2R \lambda(A) +2R \int_{A^c} |\nabla \psi-e|d\lambda \les R \int_{\R^n} |\nabla \psi-e|d\lambda\\
        &\stackrel{\eqref{firststep_stab_p1}}{\les} \eps^{-n} \delta_\eps^{\frac{1}{\alpha}}(\lambda,\mu).
        \end{aligned}
    \]
    This proves \eqref{quantmonot}.\\
We finally derive \eqref{stabmarg} from \eqref{quantmonot}. For every $1-$Lipschitz function $\vhi$ not depending on $\ang{x,e}$ we have
\[
 \int_{e^\perp} \vhi d(p_e\#\lambda)- \int_{e^\perp} \vhi d(p_e\#\mu)=\int_{\R^n} \vhi d\lambda- \int_{\R^n} \vhi d\mu \stackrel{\eqref{quantmonot}}{\les} \eps^{-n} \delta_\eps^{\frac{1}{\alpha}}(\lambda,\mu).
\]
Taking the supremum with respect to such functions  and recalling \eqref{defW1dual} we obtain \eqref{stabmarg}.

\end{proof}

\section*{Acknowledgments}
We warmly thank B. Klartag for suggesting the problem, for valuable exchanges and feedback. We thank A. Chambolle for a useful discussion regarding \eqref{estimsumxx'conv}. M.G. acknowledge partial support from the ANR Stoiques. M.F. received support under
the program "Investissement d'Avenir" launched by the French Government
and implemented by ANR, with the reference ANR‐18‐IdEx‐0001
as part of its program "Emergence".  M.F. was also supported by the Agence Nationale de la Recherche (ANR) Grant ANR-23-CE40-0003 (Project CONVIVIALITY). Part of this work was done while M.F. was a visitor at the Institute of Science and Technology Austria and the Erwin Schr\"odinger Institute in the Fall 2025. This research benefited from the support of the FMJH Program Gaspard Monge for optimization and operations research and their interactions with data science. Part of this work appears in the M.Sc. thesis of D. Tsodyks, written under the supervision of B. Klartag at the Weizmann Institute of Science.

\bibliographystyle{acm}
\bibliography{OT}

@article{FordThesis,
  title={Quantitative Stability in Discrete Optimal Transport
},
  author={Ford, Will},
  journal={arXiv preprint arXiv:2510.17407},
  year={2025}
}

@article{bolley2014dimensional,
  title={Dimensional contraction via {M}arkov transportation distance},
  author={Bolley, Fran{\c{c}}ois and Gentil, Ivan and Guillin, Arnaud},
  journal={Journal of the London Mathematical Society},
  volume={90},
  number={1},
  pages={309--332},
  year={2014},
  publisher={Oxford University Press}
}

@book{AmbDanc,
  title={Calculus of variations and partial differential equations: topics on geometrical evolution problems and degree theory},
  author={Ambrosio, Luigi and Dancer, Norman},
  year={2000},
  publisher={Springer Science \& Business Media}
}

@article{chen2022asymptotics,
  title={Asymptotics of smoothed {W}asserstein distances},
  author={Chen, Hong-Bin and Niles-Weed, Jonathan},
  journal={Potential Analysis},
  volume={56},
  number={4},
  pages={571--595},
  year={2022},
  publisher={Springer}
}

@article{BaMan,
  title={Stability Bounds for Smooth Optimal Transport Maps and their Statistical Implications},
  author={Balakrishnan, Sivaraman and Manole, Tudor},
  journal={arXiv preprint arXiv:2502.12326},
  year={2025}
}

@phdthesis{delalande2022quantitative,
  title={Quantitative Stability in Quadratic Optimal Transport},
  author={Delalande, Alex},
  year={2022},
  school={Universit{\'e} Paris-Saclay}
}

@article{DelalMeri,
  title={Quantitative stability of optimal transport maps under variations of the target measure},
  author={Delalande, Alex and Merigot, Quentin},
  journal={Duke Mathematical Journal},
  volume={172},
  number={17},
  pages={3321--3357},
  year={2023},
  publisher={Duke University Press}
}

@article{carlier2024quantitative,
  title={Quantitative stability of barycenters in the {W}asserstein space},
  author={Carlier, Guillaume and Delalande, Alex and Merigot, Quentin},
  journal={Probability Theory and Related Fields},
  volume={188},
  number={3},
  pages={1257--1286},
  year={2024},
  publisher={Springer}
}

@article{MisTre,
       author = {{Mischler}, Octave and {Trevisan}, Dario},
        title = "{Quantitative stability in optimal transport for general power costs}",
      journal = {arXiv e-prints},
         year = 2024,
        month = jul,
          eid = {arXiv:2407.19337},
        pages = {arXiv:2407.19337},
          doi = {10.48550/arXiv.2407.19337},
archivePrefix = {arXiv},
       eprint = {2407.19337},
 primaryClass = {math.FA},
       adsurl = {https://ui.adsabs.harvard.edu/abs/2024arXiv240719337M}
}

@article{kitagawa2025stability,
  title={Stability of optimal transport maps on {R}iemannian manifolds},
  author={Kitagawa, Jun and Letrouit, Cyril and M{\'e}rigot, Quentin},
  journal={arXiv preprint arXiv:2504.05412},
  year={2025}
}

@book{hajlasz2000sobolev,
  title={Sobolev met {P}oincar{\'e}},
  author={Haj{\l}asz, Piotr and Koskela, Pekka},
  volume={688},
  year={2000},
  publisher={American Mathematical Soc.}
}

@inproceedings{grigor2005stability,
  title={Stability results for {H}arnack inequalities},
  author={Grigor'yan, Alexander and Saloff-Coste, Laurent},
  booktitle={Annales de l'institut Fourier},
  volume={55},
  number={3},
  pages={825--890},
  year={2005}
}

@article{LetMer,
  title={Gluing methods for quantitative stability of optimal transport maps},
  author={Letrouit, Cyril and M{\'e}rigot, Quentin},
  journal={arXiv preprint arXiv:2411.04908},
  year={2024}
}

@inproceedings{merigot2020quantitative,
  title={Quantitative stability of optimal transport maps and linearization of the 2-{W}asserstein space},
  author={M{\'e}rigot, Quentin and Delalande, Alex and Chazal, Frederic},
  booktitle={International Conference on Artificial Intelligence and Statistics},
  pages={3186--3196},
  year={2020},
  organization={PMLR}
}

@article{brascamp1976extensions,
  title={On extensions of the {B}runn-{M}inkowski and {P}r{\'e}kopa-{L}eindler theorems, including inequalities for log concave functions, and with an application to the diffusion equation},
  author={Brascamp, Herm Jan and Lieb, Elliott H},
  journal={Journal of functional analysis},
  volume={22},
  number={4},
  pages={366--389},
  year={1976},
  publisher={Elsevier}
}

@article{HuGoTre,
  title={Asymptotics for Random Quadratic Transportation Costs},
  author={Goldman, Michael and Huesmann, Martin and Trevisan, Dario},
  journal={arXiv preprint arXiv:2409.08612},
  year={2024}
}

@book{ambrosio2005gradient,
  title={Gradient flows: in metric spaces and in the space of probability measures},
  author={Ambrosio, Luigi and Gigli, Nicola and Savar{\'e}, Giuseppe},
  year={2005},
  publisher={Springer Science \& Business Media}
}

@article{gigli2011holder,
  title={On {H}{\"o}lder continuity-in-time of the optimal transport map towards measures along a curve},
  author={Gigli, Nicola},
  journal={Proceedings of the Edinburgh Mathematical Society},
  volume={54},
  number={2},
  pages={401--409},
  year={2011},
  publisher={Cambridge University Press}
}

@article{peyre2019computational,
  title={Computational Optimal Transport: With Applications to Data Science},
  author={Peyr{\'e}, G. and Cuturi, M.},
  journal={Foundations and Trends{\textregistered} in Machine Learning},
  volume={11},
  number={5-6},
  pages={355--607},
  year={2019},
  publisher={Now Publishers, Inc.}
}

@article{de2002regularity,
  title={Regularity properties for {M}onge transport density and for solutions of some shape optimization problem},
  author={De Pascale, Luigi and Pratelli, Aldo},
  journal={Calculus of Variations and Partial Differential Equations},
  volume={14},
  number={3},
  pages={249--274},
  year={2002},
  publisher={Springer}
}

@incollection{AmbPra,
  title={Existence and stability results in the ${L}^1$ theory of optimal transportation},
  author={Ambrosio, Luigi and Pratelli, Aldo},
  booktitle={Optimal Transportation and Applications: Lectures given at the CIME Summer School, held in Martina Franca, Italy, September 2-8, 2001},
  pages={123--160},
  year={2004},
  publisher={Springer}
}

@article{hutter2021minimax,
  title={Minimax estimation of smooth optimal transport maps},
  author={H{\"u}tter, Jan-Christian and Rigollet, Philippe},
  year={2021}
}

@article{gozlan2025global,
  title={Global Regularity Estimates for Optimal Transport via Entropic Regularisation},
  author={Gozlan, Nathael and Sylvestre, Maxime},
  journal={arXiv preprint arXiv:2501.11382},
  year={2025}
}

@article{caffarelli2000monotonicity,
  title={Monotonicity properties of optimal transportation and the {FKG} and related inequalities},
  author={Caffarelli, Luis A},
  journal={Communications in Mathematical Physics},
  volume={214},
  pages={547--563},
  year={2000},
  publisher={Springer}
}

@book {AFP,
    AUTHOR = {Ambrosio, L. and Fusco, N. and Pallara, D.},
     TITLE = {Functions of bounded variation and free discontinuity
              problems},
    SERIES = {Oxford Mathematical Monographs},
 PUBLISHER = {The Clarendon Press, Oxford University Press, New York},
      YEAR = {2000},
     PAGES = {xviii+434},
      ISBN = {0-19-850245-1},
}

@article{von2005transport,
  title={Transport inequalities, gradient estimates, entropy and {R}icci curvature},
  author={von Renesse, Max-K and Sturm, Karl-Theodor},
  journal={Communications on pure and applied mathematics},
  volume={58},
  number={7},
  pages={923--940},
  year={2005},
  publisher={Wiley Online Library}
}

@article{cosenza2024new,
  title={New dimensional bounds for a branched transport problem},
  author={Cosenza, Alessandro and Goldman, Michael and Koser, Melanie},
  journal={arXiv preprint arXiv:2411.14547},
  year={2024}
}

@article{AmStTr16,
  title={A {PDE} approach to a 2-dimensional matching problem},
  author={Ambrosio, L. and Stra, F. and Trevisan, D.},
  FJournal={Probabability Theory and Related Fields},
 Journal = {{Probab. Theory Relat. Fields}}, 
volume={173},
  number={1-2},
  pages={433--477},
  year={2019},
  publisher={Springer}
}

@book {Viltop,
    AUTHOR = {Villani, C.},
     TITLE = {Topics in optimal transportation},
    SERIES = {Graduate Studies in Mathematics},
    VOLUME = {58},
 PUBLISHER = {American Mathematical Society, Providence, RI},
      YEAR = {2003},
     PAGES = {xvi+370},
      ISBN = {0-8218-3312-X},
       DOI = {10.1007/b12016},
       URL = {http://dx.doi.org/10.1007/b12016},
}

@book {Santam,
    AUTHOR = {Santambrogio, F.},
     TITLE = {Optimal transport for applied mathematicians},
    SERIES = {Progress in Nonlinear Differential Equations and their
              Applications},
    VOLUME = {87},
      NOTE = {Calculus of variations, PDEs, and modeling},
 PUBLISHER = {Birkh\"auser/Springer, Cham},
      YEAR = {2015},
     PAGES = {xxvii+353},
      ISBN = {978-3-319-20827-5; 978-3-319-20828-2},
       DOI = {10.1007/978-3-319-20828-2},
       URL = {http://dx.doi.org/10.1007/978-3-319-20828-2},
}

@article{goldman2023non,
  title={Non-convex functionals penalizing simultaneous oscillations along two independent directions: structure of the defect measure},
  author={Goldman, Michael and Merlet, Beno{\^\i}t},
  journal={arXiv preprint arXiv:2309.17067},
  year={2023}
}

@article{Van2014,
  title={R{\'e}nyi divergence and {K}ullback-{L}eibler divergence},
  author={Van Erven, Tim and Harremos, Peter},
  journal={IEEE Transactions on Information Theory},
  volume={60},
  number={7},
  pages={3797--3820},
  year={2014},
  publisher={IEEE}
}

@article{Ponce04,
  title={An estimate in the spirit of {P}oincar{\'e}'s inequality},
  author={Ponce, Augusto C},
  journal={Journal of the European Mathematical Society},
  volume={6},
  number={1},
  pages={1--15},
  year={2004}
}

@article{Bourgain01,
  title={Another look at {S}obolev spaces},
  author={Bourgain, Jean and Brezis, Haim and Mironescu, Petru},
  year={2001}
}

@article{brezis2002recognize,
  title={How to recognize constant functions. {C}onnections with {S}obolev spaces},
  author={Brezis, Ha{\i}m},
  journal={Russian Mathematical Surveys},
  volume={57},
  number={4},
  pages={693--708},
  year={2002},
  publisher={London: London Mathematical Society; distributed by Cleaver-Hume Press,[1960-}
}

@article{ollivier2009,
 author = {Ollivier, Yann},
 title = {Ricci curvature of {Markov} chains on metric spaces},
 fjournal = {Journal of Functional Analysis},
 journal = {J. Funct. Anal.},
 issn = {0022-1236},
 volume = {256},
 number = {3},
 pages = {810--864},
 year = {2009},
 language = {English},
 doi = {10.1016/j.jfa.2008.11.001},
 zbMATH = {5503526},
 Zbl = {1181.53015}
}

@article{bobkov2025,
title = {Weighted {CKP} Inequalities Involving {R}\'enyi Divergence Powers}, 
author = {Bobkov, Sergey G. and Duggal, Devraj},
journal = {Arxiv preprint},
year = {2025},
}
\end{document}